\documentclass[12pt]{article}
\usepackage{amsmath}
\usepackage{amssymb}
\usepackage{vatola}

\allowdisplaybreaks[4]
\def\q{\quad}
\def\qq{\qquad}
\def\qtq#1{\q\t{#1}\q}
\def\mod#1{\ (\text{\rm mod}\ #1)}
\def\t{\text}
\def\f{\frac}
\def\e{\equiv}
\def\b{\binom}
\def\a{\alpha}
\def\be{\beta}
\def\ga{\gamma}
\def\de{\delta}
\def\al{\alpha}
\def\ord{\t{\rm ord}}
\def\Lra{\Leftrightarrow}
\def\v2{\vskip0.2cm}

\def\sls#1#2{(\f{#1}{#2})}
 \def\ls#1#2{\big(\f{#1}{#2}\big)}
\def\Ls#1#2{\Big(\f{#1}{#2}\Big)}
\def\o{\omega}
\def\cs#1#2{\big(\f{#1}{#2}\big)_3}
\def\Cs#1#2{\Big(\f{#1}{#2}\Big)_3}
\TagsOnRight

\begin{document}

 \centerline {\bf
Cubic congruences and binary quadratic forms}
\par\q\newline \centerline{Zhi-Hong Sun}
\newline \centerline{School of Mathematics
and Statistics}\centerline{Huaiyin Normal University}
\centerline{Huaian, Jiangsu 223300, P.R. China} \centerline{Email:
zhsun@hytc.edu.cn} \centerline{URL:
http://maths.hytc.edu.cn/szh1.htm}
\par\q
\par {\it Abstract.} Let $p>3$ be a prime, $a_1,a_2,a_3\in\Bbb Z$ and let $N_p(x^3+a_1x^2+a_2x+a_3)$ denote the number of solutions to the congruence $x^3+a_1x^2+a_2x+a_3\equiv 0\pmod p$. In this paper, we give an explicit criterion for $N_p(x^3+a_1x^2+a_2x+a_3)=3$ via binary quadratic forms.


\section*{1. Introduction}

\par  For $a_1,a_2,a_3\in\Bbb Z$ let
$$\aligned&P_0=-2a_1^3+9a_1a_2-27a_3,\q Q_0=(a_1^2-3a_2)^3,
\\&D_0=-\f 1{27}(P_0^2-4Q_0)=a_1^2a_2^2-4a_2^3-4a_1^3a_3-27a_3^2+18a_1a_2a_3.
 \endaligned\tag 1.1$$
 It is known that $D_0$ is the discriminant of $x^3+a_1x^2+a_2x+a_3$.
 Let $p>3$ be a prime.
 The number of solutions to the congruence
$x^3+a_1x^2+a_2x+a_3\e 0\mod p$ is given by
$$N_p(x^3+a_1x^2+a_2x+a_3)=\begin{cases} 0\ \t{or}\ 3&\t{if $\ls {D_0}p=1,$}
\\3&\t{if $\ls {D_0}p=0$,}
\\1&\t{if $\ls {D_0}p=-1$,}\end{cases}\eqno{(1.2)}$$
where $(\f ap)$ is the Legendre symbol. See [3,9,13].

\par For $a,b,c\in\Bbb Z$ denote the binary quadratic form $ax^2+bxy+cy^2$
by $(a,b,c)$. The discriminant of $(a,b,c)$ is the number
$D=b^2-4ac$. It is well known that equivalent forms have the same discriminant and that for given form $(a,b,c)$ and positive integer $M$, there exists an equivalent form $(a',b',c')$
with $(a',M)=1$. See for example [2,16]. If $\t{gcd}(a,b,c)=1$, the form $(a,b,c)$ is said to be primitive.
Denote the equivalence class containing the form
$(a,b,c)$ by $[a,b,c]$.
If a positive integer $n$ is represented
by $(a,b,c)$, it is well known that $n$ can be represented by any form in $[a,b,c]$. Thus we may say that $n$ is represented by $[a,b,c]$. For
any nonsquare integer $D\e 0,1\mod 4$ let $H(D)$ be the form class
group consisting of classes of primitive, integral binary
quadratic forms of discriminant $D$, and let $h(D)$ be the
 class number given by $h(D)=|H(D)|$.
\par  In 1992, using class field theory Spearman and Williams [10,11]
 obtained the following result:
\par{\bf Theorem 1.1} {\it Let $a_1,a_2,a_3\in\Bbb Z$ be such that
$f(x)=x^3+a_1x^2+a_2x+a_3$ is irreducible over the field of rational numbers. Let $D_0$
be the discriminant of $f(x)$ given in (1.1). Suppose that $D_0$ is not a square and $p>3$ is a
prime with $\ls {D_0}p=1$. Then there is a unique subgroup
$J(a_1,a_2,a_3)$ of index $3$ in $H(D_0)$ such that  the congruence $f(x)\e 0\mod p$ has
three solutions if and only if $p$ is represented by one of the
classes in $J(a_1,a_2,a_3)$.}
\par In 2007, using cubic reciprocity law the author[14] proved the following theorem.
\par{\bf Theorem 1.2} {\it Let $p>3$ be a
prime, $a_1,a_2,a_3\in\Bbb Z$, $p\nmid P_0Q_0(P_0^2-4Q_0)$ and $P_0^2-4Q_0=df^2\ (d,f\in\Bbb Z)$, where $P_0$ and $Q_0$ are given by (1.1). Then the
congruence $x^3+a_1x^2+a_2x+a_3\e 0\mod p$ has three solutions if
and only if $p$ is represented by some class in $M(P_0,Q_0,f)$, where
$M(P_0,Q_0,f)$ is a subgroup of index $3$ in $H(-3k^2d)$ given in [13,Theorem 6.2]}.
\par In Theorem 1.1, $J(a_1,a_2,a_3)$ is not constructible, and in Theorem 1.2 the discriminant $-3k^2d$ is somewhat larger
than expected in some cases. In this paper, by calculating certain cubic Jacobi symbols we obtain the final criterion for  $N_p(x^3+a_1x^2+a_2x+a_3)=3$ via binary quadratic forms.
\par Let $P,D\in\Bbb Z$ with $3\nmid P$ or $27\mid P$. Suppose $Q=P^2+27D=2^{q_1}q_0^3$ with $q_0\in\Bbb Z$ and $q_1\in\{0,1,2\}$. Assume that $DQ\not=0$. For $[a,2b,c]\in H(4D)$ with $(a,3DQ)=1$ define
$\chi([a,2b,c])=\cs{P-3b(1+2\o)}a$, where $(m,n)$ is the greatest common divisor of $m$ and $n$.
In Section 3 we show that $\chi$ is
 a character on $H(4D)$ and so
  $$G(P,D)=\Big\{[a,2b,c]\in H(4D)\Bigm| (a,3DQ)=1, \Cs{P-3b(1+2\o)}a=1\Big\}$$
is a subgroup in $H(4D)$, where $a>0$ when $D<0$, $\o=(-1+\sqrt{-3})/2$ and $\cs{\alpha}{\pi}$ is the cubic Jacobi symbol defined in [7,12,14].
   Suppose that $p$ is a prime of the form $3k+1$, $p\nmid Q$ and $\ls {P^2-Q}p=1$. It is proved in Section 5 that $Q(P+\sqrt{P^2-Q})$ is a cubic residue $\mod p$ if and only if $p$ is represented by some class in $G(P,D)$.
\par   For two real numbers $b$ and $c$, the Lucas sequences $U_n=U_n(b,c)$ and $V_n=V_n(b,c)$ are defined by
 \begin{align*} &U_0=0,\ U_1=1,\ U_{n+1}=bU_n-cU_{n-1}\ (n\ge 1);
 \\&V_0=2,\ V_1=b,\ V_{n+1}=bV_n-cV_{n-1}\ (n\ge 1).
 \end{align*}
 In 2003, the author[13] introduced the third-order recurrence sequences $u_n=u_n(a_1,a_2,a_3)$ and $s_n=s_n(a_1,a_2,a_3)$ as below:
\begin{align*}&u_{-2}=u_{-1}=0,\ u_0=1,\ u_{n+3}+a_1u_{n+2}+a_2u_{n+1}+a_3u_n=0\ (n\ge -2);
\\&s_0=3,\ s_1=-a_1,\ s_2=a_1^2-2a_2,
\ s_{n+3}+a_1s_{n+2}+a_2s_{n+1}+a_3s_n=0\ (n\ge 0).\end{align*}
It is clear that $u_1=-a_1,\ u_2=a_1^2-a_2,\ u_3=-a_1^3+2a_1a_2-a_3.$
 \par For $a_1,a_2,a_3\in\Bbb Z$ let $P_0,Q_0,D_0$ be given by (1.1), and let
 $$D_1=\begin{cases} D_0/4&\t{if $2\mid P_0$,}
\\D_0&\t{if $2\nmid P_0$,}\end{cases}\  P_1=\begin{cases} P_0/2&\t{if $2\mid P_0$,}\\P_0&\t{if $2\nmid P_0$}\end{cases}\ \t{and}\ Q_1=\begin{cases} Q_0&\t{if $2\mid P_0$,}
\\4Q_0&\t{if $2\nmid P_0$.}\end{cases}\eqno{(1.3)}$$
Then $P_1^2+27D_1=Q_1$.
Let $[\alpha]$ be the greatest integer not exceeding $\alpha$.
In Section 4, we prove the following general result.
\par {\bf Theorem 1.3} Suppose that $a_1,a_2,a_3\in\Bbb Z$,
$P_0,Q_0,D_0,P_1,Q_1,D_1$ are given in (1.1) and (1.3),
 and $p>3$ is a prime such that $p\nmid D_0P_0Q_0$. Then the following statements are equivalent:
\par (i) $x^3+a_1x^2+a_2x+a_3\e 0\mod p$ has three solutions;
\par (ii) $p$ is represented by some class in $G(P_1,D_1)$;
\par (iii) $p\mid u_{p-2}(a_1,a_2,a_3)$;
\par (iv) $p\mid u_p(0,-3(a_1^2-3a_2), 2a_1^3-9a_1a_2+27a_3)$;
\par (v) $s_{p+1}(a_1,a_2,a_3)\e a_1^2-2a_2\mod p$;
\par (vi) $p\mid U_{\f{p-\sls p3}3}(P_0,Q_0)$;
\par (vii) $V_{\f{p-\sls p3}3}(P_0,Q_0)\e 2(a_1^2-3a_2)^{\f{1-\sls p3}2}
\mod p$;
\newline Moreover, if $p\nmid P_0^2-3Q_0$, then the above statements are equivalent to
\par (viii) $\sum_{k=1}^{[p/3]}\b{3k}k\sls{P_0^2}{27Q_0}^k\e 0\mod p$;
\newline if $a_1=0$, the above statements are equivalent to
\par (ix) $p\mid u_p(0,a_2,a_3)$.
\v2
 \par For any prime $p>3$, Evink and Helminck[5]
proved that $$p\mid u_{p-2}(-1,-1,-1)\iff p=x^2+11y^2\ (x,y\in\Bbb Z),$$ for $p\not=11,19$, and Faisant[6] showed that
 $$p\mid u_{p}(0,-1,-1)\iff p=x^2+23y^2\ (x,y\in\Bbb Z).$$
 Using class field theory, Moree and Noubissie[8] extended such results
 for $16$ values of $\langle a_1,a_2,a_3\rangle$.
 Clearly, our Theorem 1.3 is a vast generalization of the above special results due to Evink, Helminck, Faisant, Moree and Noubissie.

 \section*{2. The computation of $\cs{Py+3(ax+by)(1+2\o)}{ax^2+2bxy+cy^2}$}
\par\q\ Let $\Bbb Z$ be the set of integers, $\omega=(-1+\sqrt{-3})/2$ and $\Bbb Z[\o]=\{a+b\o\mid a,b\in\Bbb Z\}$.
  For $a+b\o\in\Bbb Z[\o]$
    the norm of $a+b\o$ is
given by $N(a+b\o)=(a+b\o)(a+b\o^2)=a^2-ab+b^2$. If $3\mid a-2$ and $3\mid b$, we say that $a+b\o\e 2\mod 3$ and $a+b\o$
 is primary.
\par For $a+b\o\in\Bbb Z$ with $N(a+b\o)>1$ and $a+b\o\e \pm 2\mod 3$, we
may write $a+b\o=\pm (a_1+b_1\o)\cdots(a_r+b_r\o)$, where $a_1+b_1\o,\ldots,a_r+b_r\o$
are primary primes in $\Bbb Z[\o]$. For $c+d\o\in \Bbb Z[\o]$, we can define the cubic
Jacobi symbol
$$\Cs{c+d\o}{a+b\o}=\Cs{c+d\o}{a_1+b_1\o}\cdots\Cs{c+d\o}{a_r+b_r\o},$$
where $\big(\f{c+d\o}{a_j+b_j\o})_3$ is the cubic residue character of $c+d\o$
modulo $a_j+b_j\o$ defined by
$$\Cs{c+d\o}{a_j+b_j\o}=\begin{cases} 0&\t{if $a_j+b_j\o\mid c+d\o$,}\\
\o^i&\t{if $(c+d\o)^{(N(a_j+b_j\o)-1)/3}\e \o^i\mod{a_j+b_j\o}$.}\end{cases}$$ For
 convenience we also define $\big(\f{c+d\o}1\big)_3=\big(\f{c+d\o}{-1}\big)_3=1$.
\par  According to [7, pp.\;112-115, 135,\;313] and [12,13] the
cubic Jacobi symbol has the following
properties:
\par  (2.1)  If $a,b\in\Bbb Z$ and $a+b\o\e 2\mod 3$, then
$$\Cs{\o}{a+b\o}=\o^{\f{a+b+1}3},\q
\Cs{1-\o}{a+b\o}=\o^{\f{2(a+1)}3}\qtq{and so}\Cs{1+2\o}{a+b\o}=\o^{\f b3}.$$
\par  (2.2)  If $a+b\o,c+d\o\in\Bbb Z[\o]$, $3\nmid ac$, $3\mid b$ and $3\mid d$, then
$$\Cs{c+d\o}{a+b\o}=\Cs{a+b\o}{c+d\o}.$$
  This is called general cubic reciprocity law due to Eisenstein.
\par $(\t{\rm 2.3})$ If $m,n\in\Bbb Z,\ 3\nmid m$ and $(m,n)=1$, then
$(\f nm)_3=1$.
\par $(\t{\rm 2.4})$ If $a,c,d\in\Bbb Z$ with $3\nmid a$ and $(a,N(c+d\o))=1$, then $\cs{c+d\o}a^{-1}=\cs{c+d\o^2}a$.
\par $(\t{2.5})$ If $a,b,c_1,d_1,c_2,d_2\in\Bbb Z$, $3\nmid a$ and $3\mid b$, then
$$\Cs{(c_1+d_1\o)(c_2+d_2\o)}{a+b\o}=\Cs{c_1+d_1\o}{a+b\o}\Cs{c_2+d_2\o}
{a+b\o}.$$
\par $(\t{2.6})$ If $a_1,b_1,a_2,b_2,c,d\in\Bbb Z$, $3\nmid a_1a_2$ and $b_2\e b_2\e 0\mod 3$, then
$$\Cs{c+d\o}{(a_1+b_1\o)(a_2+b_2\o)}=\Cs{c+d\o}{a_1+b_1\o}
\Cs{c+d\o}{a_2+b_2\o}.$$

\par{\bf Theorem 2.1} {\sl Let $P,D\in\Bbb Z$ with $3\nmid P$ or $27\mid P$. Suppose $Q=P^2+27D=2^{q_1}q_0^3$ with $q_0\in\Bbb Z$ and $q_1\in\{0,1,2\}$, $a,b,c,x,y\in\Bbb Z$, $(a(ax^2+2bxy+cy^2),3DQy)=1$ and $b^2-ac=D$. Then
$$\Cs{Py+3(ax+by)(1+2\o)}{ax^2+2bxy+cy^2}
=\Cs{P-3b(1+2\o)}a.$$}
\par{\it Proof.} Set $p=ax^2+2bxy+cy^2$. Then
$ap=(ax+by)^2-Dy^2$.
Since $(ap,Dy)=1$ we see that $(ap,ax+by)=1$.
Clearly
\begin{align*}&N(Py+3(ax+by)(1+2\o))
\\&=P^2y^2+27(ax+by)^2=P^2y^2+27(ap+Dy^2)=27ap+Qy^2.
\end{align*}
By the assumption, $(ap,3Qy)=1$. Thus, $(ap,N(Py+3(ax+by)(1+2\o)))=1$ and so $\cs{Py+3(ax+by)(1+2\o)}{ap}\not=0$. In addition,
$$\Cs{9b+P(1+2\o)}a
=\Cs{1+2\o}a\Cs{P-3b(1+2\o)}a=\Cs{P-3b(1+2\o)}a.$$
\par {\bf Case 1}. $27\mid Py$ and
$3\nmid ax+by$. In this case,
\begin{align*}&\Cs{Py+3(ax+by)(1+2\o)}p
\\&=\Cs{1+2\o}p\Cs{ax+by-\f{Py}9(1+2\o)}p
\\&=\Cs{ax+by-\f{Py}9-\f{2Py}9\o}{ap}
\Cs{ax+by-\f{Py}9-\f{2Py}9\o}a^{-1}
\\&=\Cs{ap}{ax+by-\f{Py}9-\f{2Py}9\o}
\Cs{by-\f{Py}9-\f{2Py}9\o}a^{-1}
\\&=\Cs{(ax+by)^2-Dy^2}{ax+by-\f{Py}9-\f{2Py}9\o}
\Cs{9b-P-2P\o}a^{-1}.
\end{align*}
Since
$$\Big(ax+by-\f{Py}9\Big)^2+\Big(ax+by-\f{Py}9\Big)\f{2Py}9
+\Ls{2Py}9^2=(ax+by)^2+\f{P^2y^2}{27},$$
for $y=3^ty_0(3\nmid y_0)$ we have
\begin{align*}&\Cs{Py+3(ax+by)(1+2\o)}p\Cs{9b-P-2P\o}a
\\&=\Cs{(ax+by)^2-Dy^2}{ax+by-\f{Py}9-\f{2Py}9\o}
=\Cs{-\f{P^2y^2}{27}-Dy^2}{ax+by-\f{Py}9-\f{2Py}9\o}
\\&=\Cs{-\f{Q}{27}\cdot 3^{2t}y_0^2}{ax+by-\f{Py}9-\f{2Py}9\o}
=\Cs {(-3)^{2t}\cdot 2^{q_1}q_0^3y_0^2}{ax+by-\f{Py}9-\f{2Py}9\o}
\\&=\Cs{\o(1-\o)}{ax+by-\f{Py}9-\f{2Py}9\o}^{4t}\Cs {ax+by-\f{Py}9-\f{2Py}9\o}2^{q_1}
\\&\q\times\Cs
{ax+by-\f{Py}9-\f{2Py}9\o}{y_0}^2
\\&=\o^{\sls{ax}3\f{2Py}{27}\cdot 4t}\Cs{ax+by-Py/9}2^{q_1}
\Cs{ax}{y_0}^2=\o^{\sls{ax}3\f{2Py}{27}\cdot 4t}=\o^{-\sls{ax}3\f{Pyt}{27}}.
\end{align*}
Thus,
\begin{align*}&\Cs{Py+3(ax+by)(1+2\o)}p
\\&=\o^{-\sls{ax}3\f{Pyt}{27}}\Cs{9b-P-2P\o}a^{-1}
=\o^{-\sls{ax}3\f{Pyt}{27}}\Cs{9b+P(1+2\o)}a
\\&=\o^{-\sls{ax}3\f{Pyt}{27}}\Cs{P-3b(1+2\o)}a.
\end{align*}
For $27\mid y$ or $27\mid P$ and $3\nmid ax+by$, we must have
$\o^{\f{Pyt}{27}}=1$ and so $\cs{Py+3(ax+by)(1+2\o)}p=\cs{P-3b(1+2\o)}a$.
\v2
\par {\bf Case 2}. $3\nmid Py$. In this case,
\begin{align*}&\Cs{Py+3(ax+by)(1+2\o)}p
\\&=\Cs{Py+3(ax+by)(1+2\o)}{ap}
\Cs{Py+3(ax+by)(1+2\o)}a^{-1}
\\&=\Cs{ap}{Py+3(ax+by)(1+2\o)}\Cs{Py+3by(1+2\o)}a^{-1}
\\&=\Cs{(ax+by)^2-Dy^2}{Py+3(ax+by)(1+2\o)}\Cs{P+3b(1+2\o)}a^{-1}
\\&=\Cs{(-3(ax+by)(1+2\o))^2+27Dy^2}{Py+3(ax+by)(1+2\o)}
\Cs{P+3b+6b\o)}a^{-1}
\\&=\Cs{P^2y^2+27Dy^2}{Py+3(ax+by)(1+2\o)}\Cs{P+3b+6b\o)}a^{-1}
\\&=\Cs{Q}{Py+3(ax+by)(1+2\o)}\Cs{y}{Py+3(ax+by)(1+2\o)}^2
\Cs{P+3b+6b\o)}a^{-1}.
\end{align*}
By (2.1)-(2.3),
\begin{align*}&\Cs{y}{Py+3(ax+by)(1+2\o)}\\&=\Cs{Py+3(ax+by)(1+2\o)}y
=\Cs{3ax(1+2\o)}y=\Cs{1+2\o}y=1.
\end{align*}
  For $2\mid Q$ we have $(ax+by)^2-Dy^2=ap\e 1\mod 2$ and so $ax+by\e 1+Dy\e 1+Py\mod 2$. Thus,
\begin{align*}\Cs{2}{Py+3(ax+by)(1+2\o)}
&=\Cs {Py+3(ax+by)(1+2\o)}2
\\&=\Cs{Py+3(ax+by)}2=1.\end{align*}
Hence we always have
$$\Big(\f{Q}{Py+3(ax+by)(1+2\o)}\Big)_3
=\Cs{2^{q_1}q_0^3}{Py+3(ax+by)(1+2\o)}=1$$
and so
$$\Cs{Py+3(ax+by)(1+2\o)}p
=\Cs{P+3b+6b\o}a^{-1}
=\Cs{P-3b-6b\o}a.$$
\par {\bf Case 3}.  $3\nmid P$, $3\mid y$, $9\nmid y$ and $ax\e -Py/3\mod 3$. In this case,
 \begin{align*}&\Cs{Py+3(ax+by)(1+2\o)}p
\\&=\Cs{\o}p\Cs{(ax+by+Py/3)\o^2+2(ax+by)}p
\\&=\Cs{\o}p\Cs{ax+by-Py/3-(ax+by+Py/3)\o}p
\\&=\Cs{\o}p\Cs{ax+by-Py/3-(ax+by+Py/3)\o}{ap}
\\&\q\times\Cs{ax+by-Py/3-(ax+by+Py/3)\o}{a}^{-1}
\\&=\Cs{\o}p\Cs{ap}{ax+by-Py/3-(ax+by+Py/3)\o}
\Cs{by-Py/3-(by+Py/3)\o}{a}^{-1}
\\&=\Cs{\o}p\Cs{(ax+by)^2-Dy^2}{ax+by-Py/3-(ax+by+Py/3)\o}
\Cs{3b-P-(3b+P)\o}{a}^{-1}.
\end{align*}
Since
\begin{align*}&\big(ax+by-\f{Py}3\big)^2+\big(ax+by-\f{Py}3\big)
\big(ax+by+\f{Py}3\big)+\big(ax+by+\f{Py}3\big)^2
\\&=3(ax+by)^2+\f{P^2y^2}{9},\end{align*}
we see that
\begin{align*}&\Cs{(ax+by)^2-Dy^2}{ax+by-\f{Py}3-(ax+by+\f{Py}3)\o}
\\&=\Cs{-27(ax+by)^2+27Dy^2}{ax+by-\f{Py}3-(ax+by+\f{Py}3)\o}
=\Cs{P^2y^2+27Dy^2}{ax+by-\f{Py}3-(ax+by+\f{Py}3)\o}
\\&=\Cs{Q\cdot 9(\f y3)^2}{ax+by-\f{Py}3-(ax+by+\f{Py}3)\o}.
\end{align*}
It is clear that
\begin{align*}
\Cs{y/3}{ax+by-\f{Py}3-(ax+by+\f{Py}3)\o}^2
&=\Cs{ax+by-\f{Py}3-(ax+by+\f{Py}3)\o}{y/3}^2
\\&=\Cs{ax(1-\o)}{y/3}^2=\Cs{1-\o}{ y/3}^2=\o^{\f{1-(\f{y/3}3)\f y3}3}
\end{align*}
and
\begin{align*}
\Cs{9}{ax+by-\f{Py}3-(ax+by+\f{Py}3)\o}
&=\Cs{\o(1-\o)}{ax+by-\f{Py}3-(ax+by+\f{Py}3)\o}^4
\\&=\o^{\f 13(\f{Py/3}3)(ax+by+\f{Py}3)}.
\end{align*}
On the other hand,
\begin{align*}&\Cs{Q}{ax+by-\f{Py}3-(ax+by+\f{Py}3)\o}
\\&=\Cs{2^{q_1}q_0^3}
{ax+by-\f{Py}3-(ax+by+\f{Py}3)\o}
=\Cs 2{ax+by-\f{Py}3-(ax+by+\f{Py}3)\o}^{q_1}
\\&=\Cs {ax+by-\f{Py}3-(ax+by+\f{Py}3)\o}2^{q_1}
=\Cs{(ax+by+\f{Py}3)(1-\o)}2^{q_1}
\\&=\Cs{1-\o}2^{q_1}=\o^{2q_1}.\end{align*}
From the above we deduce that \begin{align*}\Cs{(ax+by)^2-Dy^2}{ax+by-\f{Py}3-(ax+by+\f{Py}3)\o}
&=\Cs{Q\cdot 9(\f y3)^2}{ax+by-\f{Py}3-(ax+by+\f{Py}3)\o}
\\&=\o^{2q_1}\cdot \o^{\f 13(\f{Py/3}3)(ax+by+\f{Py}3)}
\cdot \o^{\f{1-(\f{y/3}3)\f y3}3}
\end{align*}
and so
\begin{align*}&\Cs{Py+3(ax+by)(1+2\o)}p
\\&=\Cs{\o}p\Cs{(ax+by)^2-Dy^2}{ax+by-Py/3-(ax+by+Py/3)\o}
\Cs{3b-P-(3b+P)\o}{a}^{-1}
\\&=\o^{\f 13(1-(\f p3)p)}\cdot \o^{2q_1}\cdot \o^{\f 13(\f{Py/3}3)(ax+by+\f{Py}3)}
\cdot \o^{\f{1-(\f{y/3}3)\f y3}3}\Cs{3b-P-(3b+P)\o^2}{a}
\\&=\o^{\f 13(1-(\f p3)p)}\cdot \o^{2q_1}\cdot
 \o^{\sls{P}3b+\f 13(\f{Py/3}3)(ax+\f{Py}3)}\cdot \o^{\f{1-(\f{y/3}3)\f y3}3}\Cs{6b+(3b+P)\o}a
 \\&=\o^{\f 13(1-(\f p3)p)+2q_1+\sls {P}3b}\cdot \o^{\f 13(1-ax\sls{ax}3)}\cdot \o^{\f 13(P\sls{P}3-1)}\Cs{6b+(3b+P)\o}a.
\end{align*}
Since $3\nmid P$ we have $3\nmid Q$ and so $3\nmid q_0$. It is clear that
$$q_0^3\Ls{q_0}3-1=\Big(q_0\Ls{q_0}3-1\Big)^3-3q_0\Big(\Ls{q_0}3-q_0
\Big)\e 0\mod 9.$$
Also, $(-2)^{q_1}=(1-3)^{q_1}\e 1-3q_1\mod 9$ and $\ls{q_0}3=\ls{2^{q_1}}3=(-1)^{q_1}$. Thus
$$P^2\e Q=(-1)^{q_1}(-2)^{q_1}q_0^3\e (-1)^{q_1}(1-3q_1)\Ls{q_0}3=1-3q_1\mod 9.$$
Note that
$$\f{P^2-1}3=\f{P\sls{P}3-1}3\big(P\sls{P}3+1\big)\e -\f{P\sls{P}3-1}3\mod 3.$$
We then get
$$\o^{\f{P\sls{P}3-1}3}=\o^{-\f{P^2-1}3}=\o^{q_1}.\eqno{(2.7)}$$
\par
Since $p=ax^2+2bxy+cy^2\e ax^2\e a\mod 3$ and $p=ax^2+2bxy+cy^2\e ax^2+2bxy\mod 9$, we see that for $ax\e \pm Py/3\mod 3$,
\begin{align*}
&\f{1-\sls p3p}3+\f{1-ax\sls{ax}3}3
\\&\e \f{1-\sls a3(ax^2+2bxy)}3+\f{1-ax\sls{ax}3}3
=\f{1-\sls a3ax^2}3-\f{2b}a\Ls a3ax\f y3+\f{1-ax\sls{ax}3}3
\\&\e \f{1-\sls a3a(x\sls x3-1+1)^2}3\mp 2bP\Ls y3^2+\f{1-ax\sls{ax}3}3
\\&\e \f{1-\sls a3a(1+2(x\sls x3-1))}3\mp 2bP+\f{1-ax\sls{ax}3}3
\\&=\f{1-\sls a3a}3-2a\Ls a3\f{x\sls x3-1}3\mp 2bP+\f{1-a\sls a3(x\sls x3-1+1)}3
\\&\e\f{1-\sls a3a}3-2\cdot\f{x\sls x3-1}3\mp 2b\Ls{P}3+\f{1-a\sls a3}3-a\Ls a3
\f{x\sls x3-1}3
\\&\e\f{a\sls a3-1}3\pm b\Ls{P}3\mod 3.
\end{align*}
Thus, for $ax\e -Py/3\mod 3$ we have
\begin{align*}&\Cs{Py+3(ax+by)(1+2\o)}p
\\&=\o^{\f 13(1-(\f p3)p)+2q_1+\sls {P}3b}\cdot \o^{\f 13(1-ax\sls{ax}3)}\cdot \o^{\f 13(P\sls{P}3-1)}\Cs{6b+(3b+P)\o}a
\\&=\o^{\f{a\sls a3-1}3}\Cs{6b+(3b+P)\o}a=\Cs{\o^2}a\Cs{6b+(3b+P)\o}a
\\&=\Cs{P+3b+6b\o^2}a=\Cs{P-3b(1+2\o)}a.
\end{align*}
\par {\bf Case 4.} $3\nmid P$, $3\mid y$, $9\nmid y$ and $ax\e Py/3\mod 3$. In this case,
\begin{align*}&\Cs{Py+3(ax+by)(1+2\o)}p
\\&
=\Cs{\o^2}p\Cs{(ax+by+\f{Py}3)\o
+2(ax+by)\o^2}p
\\&=\Cs{\o}p^2\Cs{-2(ax+by)+(\f{Py}3-(ax+by))\o}p
\\&=\Cs{\o}p^2\Cs{2(ax+by)+(ax+by-\f{Py}3)\o}{ap}
\Cs{2(ax+by)+(ax+by-\f{Py}3)\o}{a}^{-1}
\\&=\Cs{\o}p^2\Cs{ap}{2(ax+by)+(ax+by-\f{Py}3)\o}
\Cs{2by+(by-\f{Py}3)\o}{a}^{-1}
\\&=\Cs{\o}p^2\Cs{(ax+by)^2-Dy^2}{2(ax+by)+(ax+by-\f{Py}3)\o}
\Cs{6b+(3b-P)\o}{a}^{-1}.
\end{align*}
Since
\begin{align*}&(2(ax+by))^2-2(ax+by)\big(ax+by-\f{Py}3\big)
+\big(ax+by-\f{Py}3\big)^2
\\&=3(ax+by)^2+\f{P^2y^2}{9},\end{align*}
we see that
\begin{align*}
&\Cs{(ax+by)^2-Dy^2}{2(ax+by)+(ax+by-\f{Py}3)\o}
\\&=\Cs{-27(ax+by)^2+27Dy^2}{2(ax+by)+(ax+by-\f{Py}3)\o}
=\Cs{P^2y^2+27Dy^2}{2(ax+by)+(ax+by-\f{Py}3)\o}
\\&=\Cs{9Q(\f y3)^2}{2(ax+by)+(ax+by-\f{Py}3)\o}.\end{align*}
Note that
\begin{align*}&\Cs{(y/3)^2}{2(ax+by)+(ax+by-Py/3)\o}
\\&=\Cs {2(ax+by)+(ax+by-Py/3)\o}{y/3}^2
=\Cs{ax(2+\o)}{y/3}^2=\Cs{-\o(1+2\o)}{y/3}^2
\\&=\Cs{\o}{y/3}^2=\o^{\f 23(1-(\f y3)\f y3)}
\end{align*}
and
\begin{align*}&\Cs{9}{2(ax+by)+(ax+by-\f{Py}3)\o}
\\&=\Cs{\o(1-\o)}{2(ax+by)+(ax+by-\f{Py}3)\o}^4
= \Cs{\o(1-\o)}{2(ax+by)\sls{ax}3+(ax+by-\f{Py}3)\sls{ax}3\o}
\\&=\o^{\f 13\sls{ax}3(ax+by-\f{Py}3)}=\o^{(\f{Py/3}3)b\f y3+\f 13
\sls{ax}3(ax-\f{Py}3)}=\o^{\sls{P}3b+\f 13
\sls{ax}3(ax-\f{Py}3)}.
\end{align*}
For $2\mid Q$ we see that $ax+by\e (ax+by)^2\e Dy^2+1\e Py^2+1\e \f{Py}3+1\mod 2$ and so
\begin{align*}\Cs{2}{2(ax+by)+(ax+by-\f{Py}3)\o}
&=\Cs{2(ax+by)+(ax+by-\f{Py}3)\o}2
\\&=\Cs {\o}2=\o.\end{align*}
Hence, we always have $$\Cs{Q}{2(ax+by)+(ax+by-\f{Py}3)\o}
=\Cs{2^{q_1}q_0^3}{2(ax+by)+(ax+by-\f{Py}3)\o}=\o^{q_1}.$$
From the above we deduce that
$$\Cs{(ax+by)^2-Dy^2}{2(ax+by)+(ax+by-\f{Py}3)\o}
=\o^{\sls{P}3b+\f 13
\sls{ax}3(ax-\f{Py}3)}\cdot \o^{q_1}\cdot
\o^{\f 23(1-(\f y3)\f y3)}$$
and so
\begin{align*}&\Cs{Py+3(ax+by)(1+2\o)}p
\\&=\Cs{\o}p^2\Cs{(ax+by)^2-Dy^2}{2(ax+by)+(ax+by-\f{Py}3)\o}
\Cs{6b+(3b-P)\o}{a}^{-1}
\\&=\o^{\f 23(1-\sls p3p)}\cdot \o^{\sls{P}3b+\f 13
\sls{ax}3(ax-\f{Py}3)}\cdot \o^{q_1}\cdot
\o^{\f 23(1-(\f {y/3}3)\f y3)}\Cs{6b+(3b-P)\o^2}{a}
\\&=\o^{-\f 13(1-\sls p3p)}\cdot \o^{\sls{P}3b}\cdot \o^{q_1}
\cdot\o^{\f 13(ax\sls{ax}3-\sls{Py/3}3\f{Py}3-1+\sls {y/3}3\f y3)}
\Cs{6b+(3b-P)\o^2}{a}
\\&=\o^{-\f 13(1-\sls p3p)}\cdot \o^{\sls{P}3b}\cdot \o^{q_1}
\cdot\o^{\f 13(ax\sls{ax}3-1+1-\sls{P}3P)}
\Cs{6b+(3b-P)\o^2}{a}.
\end{align*}
From the arguments in Case 3,
$$ \o^{\f{P\sls{P}3-1}3}=\o^{q_1}\q\t{and}\q
\q \o^{\f{1-\sls p3p}3+\f{1-ax\sls{ax}3}3}=\o^{\f {a\sls a3-1}3+b\sls{P}3}.$$
Therefore,
\begin{align*}&\Cs{Py+3(ax+by)(1+2\o)}p
\\&=\o^{-\f 13(1-\sls p3p)}\cdot \o^{\sls{P}3b}\cdot \o^{q_1}
\cdot\o^{\f 13(ax\sls{ax}3-1+1-\sls{P}3P)}
\Cs{6b+(3b-P)\o^2}{a}
\\&=\o^{\f {1-a\sls a3}3}\Cs{6b+(3b-P)\o^2}{a}
=\Cs {-\o}a\Cs{6b+(3b-P)\o^2}{a}
\\&=\Cs{P-3b(1+2\o)}a.
\end{align*}

\par {\bf Case 5.} $3\nmid P$, $9\mid y$, $27\nmid y$ and $ax\e \f{Py}{9}\mod 3$. In this case,
\begin{align*}&\Cs{Py+3(ax+by)(1+2\o)}p
\\&=\Cs{1+2\o}p\Cs{ax+by-\f{Py}{9}(1+2\o)}p
=\Cs{ax+by-\f{Py}{9}-\f{2Py}9\o}p
\\&=\Cs{\o}p\Cs{-\f{2Py}9+(ax+by-\f{Py}{9})\o^2}p
=\Cs{\o}p\Cs{-\f{2Py}9+(ax+by-\f{Py}{9})\o}p^{-1}
\\&=\Cs{\o}p\Cs{-\f{2Py}9+(ax+by-\f{Py}{9})\o}{ap}^{-1}
\Cs{-\f{2Py}9+(ax+by-\f{Py}{9})\o}{a}
\\&=\Cs{\o}p\Cs{ap}{-\f{2Py}9+(ax+by-\f{Py}{9})\o}^{-1}
\Cs{-\f{2Py}9+(by-\f{Py}{9})\o}{a}
\\&=\o^{\f 13(1-(\f p3)p)}\Cs{(ax+by)^2-Dy^2}{-\f{2Py}9+(ax+by-\f{Py}{9})\o}^{-1}
\Cs{-2P+(9b-P)\o}{a}.
\end{align*}
Since
$$\Big(-\f{2Py}9\Big)^2+\f{2Py}9\Big(ax+by-\f{Py}{9}\Big)
+\Big(ax+by-\f{Py}{9}\Big)^2=(ax+by)^2+\f{P^2y^2}{27},$$
we see that
\begin{align*}&\Cs{(ax+by)^2-Dy^2}{-\f{2Py}9+(ax+by-\f{Py}{9})\o}
\\&=\Cs{-\f{P^2y^2}{27}-Dy^2}{-\f{2Py}9+(ax+by-\f{Py}{9})\o}
=\Cs{-3\cdot 2^{q_1}q_0^3\sls y9^2}{-\f{2Py}9+(ax+by-\f{Py}{9})\o}
\\&=\Cs{\o(1-\o)}{-\f{2Py}9+(ax+by-\f{Py}{9})\o}^2
\Cs 2{-\f{2Py}9+(ax+by-\f{Py}{9})\o}^{q_1}
\\&\q\times\Cs{\f y9}{-\f{2Py}9+(ax+by-\f{Py}{9})\o}^2
\\&=\o^{-\f 23(\f{Py/9}3)(ax+by-\f{Py}{9})}
\Cs {-\f{2Py}9+(ax+by-\f{Py}{9})\o}2^{q_1}
\\&\q\times\Cs{-\f{2Py}9
+(ax+by-\f{Py}{9})\o}{\f y9}^2
\\&=\o^{-\f 23(\f{Py/9}3)(ax-\f{Py}{9})}\Cs{\o}2^{q_1}
\Cs{\o}{y/9}^2
=\o^{\f 13(\f{Py/9}3)(ax-\f{Py}9)}\cdot \o^{q_1}\cdot
\o^{\f 23(1-(\f{y/9}3)\f y9)}.
\end{align*}
and so
\begin{align*}&\Cs{Py+3(ax+by)(1+2\o)}p
\\&=\o^{\f 13(1-(\f p3)p)}\Cs{(ax+by)^2-Dy^2}{-\f{2Py}9+(ax+by-\f{Py}{9})\o}^{-1}
\Cs{-2P+(9b-P)\o}{a}
\\&=\o^{\f 13(1-(\f p3)p)}\cdot \o^{-\f 13(\f{Py/9}3)(ax-\f{Py}9)}\cdot \o^{-q_1}\cdot
\o^{\f 13(1-(\f{y/9}3)\f y9)}\Cs{-2P+(9b-P)\o}{a}
\\&=\o^{\f 13(1-\sls a3ax^2)}\cdot\o^{\f 13(1-ax\sls{ax}3+(\f{y/9}3)\f y9(P\sls{P}3-1))}\cdot \o^{-q_1}\Cs{-2P+(9b-P)\o}{a}
\\&=\o^{\f 13(2-ax\sls{ax}3(2+x\sls x3-1))}\cdot \o^{\f 13(P\sls{P}3-1)}\cdot\o^{-q_1}\Cs{-2P+(9b-P)\o}{a}
\\&=\o^{\f 13(ax\sls{ax}3-1+1-x\sls x3)}\cdot \o^{\f 13(P\sls{P}3-1)}\cdot\o^{-q_1}\Cs{-2P+(9b-P)\o}{a}
\\&=\o^{-\f 13(1-a\sls a3)}\cdot \o^{\f 13(P\sls{P}3-1)}\cdot\o^{-q_1}\Cs{-2P+(9b-P)\o}{a}.
\end{align*}
Since $ \o^{\f{P\sls{P}3-1}3}=\o^{q_1}$ by (2.7), we deduce that
\begin{align*}&\Cs{Py+3(ax+by)(1+2\o)}p
\\&=\o^{-\f 13(1-a\sls a3)}\Cs{-2P+(9b-P)\o}{a}
=\Cs{\o^2}a\Cs{-2P+(9b-P)\o}{a}
\\&=\Cs{9b-P-2P\o^2}a=\Cs{9b+P(1+2\o)}a
=\Cs{P-3b(1+2\o)}a.
\end{align*}
\par{\bf Case 6.} $3\nmid P$, $9\mid y$, $27\nmid y$ and $ax\e -Py/9\mod 3$. In this case,
\begin{align*}&\Cs{Py+3(ax+by)(1+2\o)}p
\\&=\Cs{1+2\o}p\Cs{ax+by-\f{Py}{9}(1+2\o)}p=
\Cs{ax+by-\f{Py}{9}-\f{2Py}9\o}p
\\&=\Cs {\o^2}p\Cs{(ax+by-\f{Py}{9})\o-\f{2Py}9(-1-\o)}p
=\Cs{\o}p^2\Cs{\f {2Py}9+(ax+by+\f{Py}{9})\o}p
\\&=\Cs{\o}p^2\Cs{\f {2Py}9+(ax+by+\f{Py}{9})\o}{ap}\Cs{\f {2Py}9+(ax+by+\f{Py}{9})\o}a^{-1}
\\&=\Cs{\o}p^2\Cs{ap}{\f {2Py}9+(ax+by+\f{Py}{9})\o}\Cs{\f {2Py}9+(by+\f{Py}{9})\o}a^{-1}
\\&=\o^{\f 23(1-(\f p3)p)}\Cs{(ax+by)^2-Dy^2}{\f {2Py}9+(ax+by+\f{Py}{9})\o}
\Cs{2P+(9b+P)\o}a^{-1}.
\end{align*}
Since
$$\Ls {2Py}9^2-\f{2Py}9\Big(ax+by+\f{Py}{9}\Big)
+\Big(ax+by+\f{Py}{9}\Big)^2
=(ax+by)^2+\f{P^2y^2}{27},$$
we see that
\begin{align*}
&\Cs{(ax+by)^2-Dy^2}{\f {2Py}9+(ax+by+\f{Py}{9})\o}
\\&=\Cs{-\f{P^2y^2}{27}-Dy^2}{\f {2Py}9+(ax+by+\f{Py}{9})\o}
=\Cs{-3\cdot 2^{q_1}q_0^3\sls y9^2}{\f {2Py}9+(ax+by+\f{Py}{9})\o}
\\&=\Cs{\o(1-\o)}{\f {2Py}9+(ax+by+\f{Py}{9})\o}^2
\Cs 2{\f {2Py}9+(ax+by+\f{Py}{9})\o}^{q_1}
\\&\q\times\Cs {\f {2Py}9+(ax+by+\f{Py}{9})\o}{\f y9}^2
\\&=\o^{\f 23(\f{Py/9}3)(ax+by+\f{Py}{9})}
\Cs {\f {2Py}9+(ax+by+\f{Py}{9})\o}2^{q_1}
\Cs {ax\o}{y/9}^2
\\&=\o^{\f 23(-ax\sls{ax}3+\sls{Py/9}3\f{Py}9)}
\Cs {\o}2^{q_1}\Cs {\o}{y/9}^2
\\&=\o^{\f 23(-ax\sls{ax}3+P\sls{P}3\f y9\sls{y/9}3)}
\cdot \o^{q_1}\cdot
\o^{\f 23(1-\f y9(\f{y/9}3))}
=\o^{\f 23(1-ax\sls{ax}3+P\sls{P}3-1)}\cdot \o^{q_1}
\end{align*}
and so
\begin{align*}&\Cs{Py+3(ax+by)(1+2\o)}p\Cs{2P+(9b+P)\o}a
\\&=\o^{\f 23(1-(\f p3)p)}\Cs{(ax+by)^2-Dy^2}{\f {2Py}9+(ax+by+\f{Py}{9})\o}
\\&=\o^{\f 23(1-(\f a3)ax^2)}\cdot\o^{\f 23(1-ax\sls{ax}3)}\cdot \o^{\f 23(P\sls{P}3-1)}\cdot \o^{q_1}.
\end{align*}
By (2.7), $\o^{\f 23(P\sls{P}3-1)}\cdot \o^{q_1}=1$. Also,
\begin{align*}\o^{\f 23(1-(\f a3)ax^2)}\cdot\o^{\f 23(1-ax\sls{ax}3)}
&=\o^{\f 23(2-2ax\sls{ax}3-ax\sls{ax}3(x\sls x3-1))}
=\o^{\f 13(1-ax\sls{ax}3+x\sls x3-1)}
\\&=\o^{\f 13x\sls x3(1-a\sls a3)}
=\o^{\f 13(1-a\sls a3)}.
\end{align*}
Thus,
\begin{align*}&\Cs{Py+3(ax+by)(1+2\o)}p
\\&=\o^{\f 13(1-a\sls a3)}
\Cs{2P+(9b+P)\o}a^{-1}=\Cs {\o}a
\Cs{2P+(9b+P)\o^2}a
\\&=\Cs{9b+P(1+2\o)}a=\Cs{P-3b(1+2\o)}a.
\end{align*}

\par {\bf Case 7}. $27\mid P$ and $3\mid ax+by$. Since $3\mid (ax+by)$ and $(ax,y)=1$ we see that $3\nmid y$.
Suppose $r=\ord_3 P$ and $s=\ord_3 (ax+by)$, where $\ord_3n$ is the great non-negative integer such that $3^{\ord_3n}\mid n$.
Since $27\mid P$ and $P^2+27D=Q=2^{q_1}q_0^3$, we have $3\mid q_0$. Set $q_2=q_0/3$. Then $P\cdot \f{P}{27}+D=2^{q_1}q_2^3$.
Note that $Dy^2=(ax+by)^2-ap\e -ap\not\e 0\mod 3$. We see that $3\nmid D$ and so $3\nmid q_2$. Thus,
\begin{align*}D&\e 2^{q_1}q_2^3=(-1)^{q_1}(1-3)^{q_1}\Ls{q_2}3\Big(1+q_2\Ls {q_2}3-1\Big)^3
\\&\e (-1)^{q_1}\Ls{q_2}3(1-3q_1)\mod 9,
\end{align*}
which implies that
$$\o^{\f{1-D\sls {D}3}3}=\o^{\f{1-(-1)^{q_1}\sls{q_2}3D}3}=\o^{q_1}.\eqno{(2.8)}$$
Observing that $ap=(ax+by)^2-Dy^2\e -Dy^2\mod 9$ we derive that
 \begin{align*}\Cs{\o}{ap}^2\cdot  \o^{\f{1-y^2}3}&=\o^{\f{2(1-ap\sls{ap}3)}3}
 \cdot  \o^{\f{1-y^2}3}=\o^{\f{2(1-Dy^2\sls{D}3)}3+\f{1-y^2}3}
 \\&=\o^{\f{D\sls{D}3-1}3y^2}=\o^{\f{D\sls{D}3-1}3}.
 \end{align*}
 Thus,
 $$\Cs{\o}{ap}^2\cdot  \o^{\f{1-y^2}3}=\o^{-q_1}
 \q\t{and so}\q\Cs{\o}{ap}\cdot\o^{-\f{1-y^2}3}=\o^{q_1}
 .\eqno{(2.9)}$$

 We first assume that $r\le s$. Observe that
$$N\Big(\f{Py}{3^r}+3\f{ax+by}{3^r}(1+2\o)\Big)=
\Ls{Py}{3^r}^2+3\Big(\f{3(ax+by)}{3^r}\Big)^2=
\f{P^2y^2+27(ax+by)^2}{3^{2r}}.$$
Using (2.1)-(2.6) we deduce that
\begin{align*}&\Cs{Py+3(ax+by)(1+2\o)}p
\\&=\Cs{\f{Py}{3^r}+3\f{ax+by}{3^r}(1+2\o)}p
=\Cs{\f{Py}{3^r}+3\f{ax+by}{3^r}(1+2\o)}{ap}
\Cs{\f{Py}{3^r}+3\f{ax+by}{3^r}(1+2\o)}a^{-1}
\\&=\Cs{ap}{\f{Py}{3^r}+3\f{ax+by}{3^r}(1+2\o)}
\Cs{Py+3(ax+by)(1+2\o)}a^{-1}
\\&=\Cs{(ax+by)^2-Dy^2}{\f{Py}{3^r}+3\f{ax+by}{3^r}(1+2\o)}
\Cs{Py+3by(1+2\o)}a^{-1}
\\&=\Cs{27(ax+by)^2-27Dy^2}{\f{Py}{3^r}+3\f{ax+by}{3^r}(1+2\o)}
\Cs{P+3b(1+2\o)}a^{-1}
\\&=\Cs {-(P^2+27D)y^2}{\f{Py}{3^r}+3\f{ax+by}{3^r}(1+2\o)}
\Cs{P-3b(1+2\o)}a
\\&=\Cs{-2^{q_1}q_0^3}{\f{Py}{3^r}+3\f{ax+by}{3^r}(1+2\o)}
\Cs {\f{Py}{3^r}+3\f{ax+by}{3^r}(1+2\o)}{y^2}\Cs{P-3b(1+2\o)}a
\\&=\Cs {2}{\f{Py}{3^r}+3\f{ax+by}{3^r}(1+2\o)}^{q_1}
\Cs {\f{Py}{3^r}+3\f{ax+by}{3^r}(1+2\o)}y^2\Cs{P-3b(1+2\o)}a
\\&=\Cs{\f{Py}{3^r}+3\f{ax+by}{3^r}(1+2\o)}2^{q_1}
\Cs{Py+3(ax+by)(1+2\o)}y^2\Cs{P-3b(1+2\o)}a
\\&=\Cs{Py+3(ax+by)}2^{q_1}\Cs{3ax(1+2\o)}y^2\Cs{P-3b(1+2\o)}a
\\&=\Cs{1+2\o}y^2\Cs{P-3b(1+2\o)}a=\Cs{P-3b(1+2\o)}a.
\end{align*}
\par Next we assume $r>s+2$. Note that
$$N\Big(\f{ax+by}{3^s}-\f{Py}{3^{s+2}}(1+2\o)\Big)
=\Ls{ax+by}{3^s}^2+3\Ls{Py}{3^{s+2}}^2
=\f{(ax+by)^2+\f{P^2y^2}{27}}{3^{2s}}.$$
Using (2.1)-(2.6) we see that
\begin{align*}&\Cs{Py+3(ax+by)(1+2\o)}p
\\&=\Cs{1+2\o}p\Cs{ax+by-\f{Py}9(1+2\o)}p
=\Cs{\f{ax+by}{3^s}-\f{Py}{3^{s+2}}(1+2\o)}p
\\&=\Cs{\f{ax+by}{3^s}-\f{Py}{3^{s+2}}(1+2\o)}{ap}
\Cs{\f{ax+by}{3^s}-\f{Py}{3^{s+2}}(1+2\o)}a^{-1}
\\&=\Cs{ap}{\f{ax+by}{3^s}-\f{Py}{3^{s+2}}(1+2\o)}
\Cs{ax+by-\f{Py}9(1+2\o)}a^{-1}
\\&=\Cs{(ax+by)^2-Dy^2}{\f{ax+by}{3^s}-\f{Py}{3^{s+2}}(1+2\o)}
\Cs{by-\f{Py}9(1+2\o)}a^{-1}
\\&=\Cs{-\f{P^2+27D}{27}y^2}{\f{ax+by}{3^s}-\f{Py}{3^{s+2}}(1+2\o)}
\Cs{9b-P(1+2\o)}a^{-1}
\\&=\Cs{2^{q_1}q_0^3y^2}{\f{ax+by}{3^s}-\f{Py}{3^{s+2}}(1+2\o)}
\Cs{9b+P(1+2\o)}a
\\&=\Cs 2{\f{ax+by}{3^s}-\f{Py}{3^{s+2}}(1+2\o)}^{q_1}
\Cs y{\f{ax+by}{3^s}-\f{Py}{3^{s+2}}(1+2\o)}^2\Cs{9b+P(1+2\o)}a
\\&=\Cs {\f{ax+by}{3^s}-\f{Py}{3^{s+2}}(1+2\o)}2^{q_1}
\Cs{\f{ax+by}{3^s}-\f{Py}{3^{s+2}}(1+2\o)}y^2\Cs{9b+P(1+2\o)}a
\\&=\Cs{ax+by-\f{P}9y(1+2\o)}y^2
\Cs{1+2\o}a\Cs{P-3b(1+2\o)}a
\\&=\Cs{ax}y^2\Cs{P-3b(1+2\o)}a=\Cs{P-3b(1+2\o)}a.
\end{align*}
\par Now suppose $r=s+2$. If $\f{ax+by}{3^s}\e\f{Py}{3^{s+2}}\mod 3$, since
$$N\Big(-\f{2Py}{3^{s+2}}+\Big(\f{ax+by}{3^s}-\f{Py}{3^{s+2}}\Big)\o\Big)
=\f {P^2y^2/27+(ax+by)^2}{3^{2s+6}},$$
using (2.1)-(2.6) and (2.9) we derive that
\begin{align*}&\Cs{Py+3(ax+by)(1+2\o)}p
\\&=\Cs{1+2\o}p\Cs{ax+by-\f{Py}9(1+2\o)}p
=\Cs{\f{ax+by}{3^s}-\f{Py}{3^{s+2}}(1+2\o)}p
\\&=\Cs{\o}p\Cs{-\f{2Py}{3^{s+2}}
+(\f{ax+by}{3^s}-\f{Py}{3^{s+2}})\o^2}p
=\Cs {\o}p\Cs{-\f{2Py}{3^{s+2}}
+(\f{ax+by}{3^s}-\f{Py}{3^{s+2}})\o}p^{-1}
\\&=\Cs {\o}p\Cs{-\f{2Py}{3^{s+2}}
+(\f{ax+by}{3^s}-\f{Py}{3^{s+2}})\o}{ap}^{-1}
\Cs{-\f{2Py}{3^{s+2}}
+(\f{ax+by}{3^s}-\f{Py}{3^{s+2}})\o}a
\\&=\Cs {\o}p\Cs{ap}{-\f{2Py}{3^{s+2}}
+(\f{ax+by}{3^s}-\f{Py}{3^{s+2}})\o}^{-1}
\Cs{-2Py+(9(ax+by)-Py)\o}a
\\&=\Cs {\o}p\Cs{(ax+by)^2-Dy^2}{-\f{2Py}{3^{s+2}}
+(\f{ax+by}{3^s}-\f{Py}{3^{s+2}})\o}^{-1}
\Cs{-2Py+(9b-P)y\o}a
\\&=\Cs {\o}p\Cs{-\f{P^2+27D}{27}y^2}{-\f{2Py}{3^{s+2}}
+(\f{ax+by}{3^s}-\f{Py}{3^{s+2}})\o}^{-1}
\Cs{2P-(9b-P)\o}a
\\&=\Cs {\o}p\Cs{2^{q_1}q_0^3y^2}{-\f{2Py}{3^{s+2}}
+(\f{ax+by}{3^s}-\f{Py}{3^{s+2}})\o}^{-1}\Cs{\o}a
\Cs{P-9b+2P\o^2}a
\\&=\Cs {\o}{ap}\Cs {-\f{2Py}{3^{s+2}}
+(\f{ax+by}{3^s}-\f{Py}{3^{s+2}})\o}2^{-q_1}\Cs {-\f{2Py}{3^{s+2}}
+(\f{ax+by}{3^s}-\f{Py}{3^{s+2}})\o}y^{-2}
\\&\q\times\Cs{P-9b+2P\o}a^{-1}
\\&=\Cs{\o}{ap}\Cs{\o}2^{-q_1}\Cs{\o}y^{-2}
\Cs{1+2\o}a^{-1}\Cs{P+3b(1+2\o)}a^{-1}
\\&=\Cs{\o}{ap}\cdot \o^{-q_1}\cdot
\o^{-\f{1-y^2}3}\Cs{P-3b(1+2\o)}a
=\Cs{P-3b(1+2\o)}a.
\end{align*}
 \par If $r=s+2$ and $\f{ax+by}{3^s}\e-\f{Py}{3^{s+2}}\mod 3$,
 using (2.1)-(2.6) and (2.9) we derive that
 \begin{align*}&\Cs{Py+3(ax+by)(1+2\o)}p
\\&=\Cs{1+2\o}p\Cs{ax+by-\f{Py}9(1+2\o)}p
\\&=\Cs{\f{ax+by}{3^s}-\f{Py}{3^{s+2}}(1+2\o)}p
=\Cs{\f{ax+by}{3^s}+\f{Py}{3^{s+2}}(1+2\o)}p^{-1}
\\&=\Cs{\o}p^{-1}\Cs{\f{2Py}{3^{s+2}}
+(\f{ax+by}{3^s}+\f{Py}{3^{s+2}})\o^2}p^{-1}
=\Cs {\o}p^2\Cs{\f{2Py}{3^{s+2}}
+(\f{ax+by}{3^s}+\f{Py}{3^{s+2}})\o}p
\\&=\Cs {\o}p^2\Cs{\f{2Py}{3^{s+2}}
+(\f{ax+by}{3^s}+\f{Py}{3^{s+2}})\o}{ap}\Cs{\f{2Py}{3^{s+2}}
+(\f{ax+by}{3^s}+\f{Py}{3^{s+2}})\o}{a}^{-1}
\\&=\Cs {\o}p^2\Cs{(ax+by)^2-Dy^2}{\f{2Py}{3^{s+2}}
+(\f{ax+by}{3^s}+\f{Py}{3^{s+2}})\o}\Cs{2Py+(9(ax+by)+Py)\o}a^{-1}
\\&=\Cs {\o}p^2\Cs{-\f{P^2}{27}y^2-Dy^2}{\f{2Py}{3^{s+2}}
+(\f{ax+by}{3^s}+\f{Py}{3^{s+2}})\o}\Cs{2Py+(9by+Py)\o}a^{-1}
\\&=\Cs {\o}p^2\Cs{2^{q_1}q_0^3y^2}{\f{2Py}{3^{s+2}}
+(\f{ax+by}{3^s}+\f{Py}{3^{s+2}})\o}\Cs{2P+(9b+P)\o}a^{-1}
\\&=\Cs {\o}p^2\Cs{\f{2Py}{3^{s+2}}
+(\f{ax+by}{3^s}+\f{Py}{3^{s+2}})\o}2^{q_1}
\Cs {\f{2Py}{3^{s+2}}
+(\f{ax+by}{3^s}+\f{Py}{3^{s+2}})\o}y^2
\\&\q\times\Cs{2P+(9b+P)\o^2}a
\\&=\Cs {\o}p^2\Cs{2Py
+(9(ax+by)+Py)\o}2^{q_1}
\Cs {2Py+(9(ax+by)+Py)\o}y^2
\\&\q\times\Cs{2P+(9b+P)\o^2}a
\\&=\Cs {\o}p^2\Cs{\o}2^{q_1}\Cs{\o}y^2
\Cs{2P+(9b+P)\o^2}a
\\&=\Cs{\o}{ap}^2\cdot \o^{q_1}\cdot \o^{\f{1-y^2}3}\Cs{P+9b+2P\o}a
=\Cs{P+9b+2P\o}a
\\&=\Cs{1+2\o}a\Cs{P-3b(1+2\o)}a=\Cs{P-3b(1+2\o)}a.
\end{align*}

 \par Finally suppose $r=s+1$. If $\f{Py}{3^{s+1}}\e -\f{ax+by}{3^s}\mod 3$, appealing to (2.1)-(2.6) and (2.9) we deduce that
 \begin{align*}&\Cs{Py+3(ax+by)(1+2\o)}p
 \\&=\Cs{\f{Py}{3^{s+1}}+\f{ax+by}{3^s}(1+2\o)}p
 =\Cs {\o}p\Cs{(\f{Py}{3^{s+1}}+\f{ax+by}{3^s})\o^2+2\f{ax+by}{3^s}}p
 \\&=\Cs{\o}p\Cs{2\f{ax+by}{3^s}+(\f{Py}{3^{s+1}}+\f{ax+by}{3^s})\o}p^{-1}
 \\&=\Cs{\o}p\Cs{2\f{ax+by}{3^s}+(\f{Py}{3^{s+1}}+\f{ax+by}{3^s})\o}
 {ap}^{-1}\Cs{2\f{ax+by}{3^s}+(\f{Py}{3^{s+1}}+\f{ax+by}{3^s})\o}
 {a}
 \\&=\Cs{\o}p\Cs{(ax+by)^2-Dy^2}
 {2\f{ax+by}{3^s}+(\f{Py}{3^{s+1}}+\f{ax+by}{3^s})\o}^{-1}
 \Cs{2(ax+by)+(\f{P}3y+ax+by)\o}
 {a}
 \\&=\Cs{\o}p\Cs{-\f{P^2+27D}{27}y^2}
 {2\f{ax+by}{3^s}+(\f{Py}{3^{s+1}}+\f{ax+by}{3^s})\o}^{-1}
 \Cs{2by+(\f{P}3y+by)\o}a
 \\&=\Cs{\o}p\Cs{2^{q_1}q_0^3y^2}
 {2\f{ax+by}{3^s}+(\f{Py}{3^{s+1}}+\f{ax+by}{3^s})\o}^{-1}
 \Cs{6b+(P+3b)\o}a
 \\&=\Cs{\o}p\Cs{2}{2\f{ax+by}{3^s}+(\f{Py}{3^{s+1}}
 +\f{ax+by}{3^s})\o}^{-q_1}
 \Cs y{2\f{ax+by}{3^s}+(\f{Py}{3^{s+1}}+\f{ax+by}{3^s})\o}^{-2}
 \\&\q\times\Cs{6b+(P+3b)\o}a
 \\&=\Cs{\o}p\Cs {2\f{ax+by}{3^s}+(\f{Py}{3^{s+1}}+\f{ax+by}{3^s})\o}2^{-q_1}
 \Cs{2\f{ax+by}{3^s}+(\f{Py}{3^{s+1}}+\f{ax+by}{3^s})\o}y^{-2}
 \\&\q\times\Cs{6b+(P+3b)\o}a
 \\&=\Cs{\o}p\Cs{\o}2^{-q_1}\Cs{2(ax+by)+(\f{P}3y+ax+by)\o}y
 \Cs{6b+(P+3b)\o}a
 \\&=\Cs{\o}p\cdot \o^{-q_1}\Cs{ax(2+\o)}y\Cs{6b+(P+3b)\o}a
 \\&=\Cs{\o}{ap}\cdot \o^{-q_1}\Cs{\o(1+2\o)}y\Cs{P+3b+6b\o^2}a
 \\&=\Cs{\o}{ap}\cdot \o^{-q_1}\cdot \o^{-\f{1-y^2}3}\Cs{P+3b+6b\o}a^{-1}
 =\Cs{P-3b(1+2\o)}a.
  \end{align*}
\par If $r=s+1$ and
 $\f{Py}{3^{s+1}}\e \f{ax+by}{3^s}\mod 3$, using (2.1)-(2.6) and (2.9) we see that
 \begin{align*}&\Cs{Py+3(ax+by)(1+2\o)}p
 \\&=\Cs{\f{Py}{3^{s+1}}+\f{ax+by}{3^s}(1+2\o)}p
 =\Cs {\o^2}p\Cs{-2\f{2x+by}{3^s}+(\f{Py}{3^{s+1}}-
 \f{ax+by}{3^s})\o}p
 \\&=\Cs{\o}p^2\Cs{-2\f{ax+by}{3^s}+(\f{Py}{3^{s+1}}-
 \f{ax+by}{3^s})\o}{ap} \Cs{-2\f{ax+by}{3^s}+(\f{Py}{3^{s+1}}-
 \f{ax+by}{3^s})\o}{a}^{-1}
 \\&=\Cs{\o}p^2
 \Cs{(ax+by)^2-Dy^2}{-2\f{ax+by}{3^s}+(\f{Py}{3^{s+1}}-
 \f{ax+by}{3^s})\o}\Cs{-2(ax+by)+(\f{P}3y-(ax+by))\o}{a}^{-1}
 \\&=\Cs{\o}p^2
 \Cs{-\f{P^2+27D}{27}y^2}{-2\f{ax+by}{3^s}+(\f{Py}{3^{s+1}}-
 \f{ax+by}{3^s})\o}\Cs{-2by+(\f{P}3y-by))\o}{a}^{-1}
 \\&=\Cs{\o}p^2
 \Cs{2^{q_1}q_0^3y^2}{-2\f{ax+by}{3^s}+(\f{Py}{3^{s+1}}-
 \f{ax+by}{3^s})\o}\Cs{-6b+(P-3b)\o}{a}^{-1}
 \\&=\Cs{\o}p^2\Cs{-2\f{ax+by}{3^s}+(\f{Py}{3^{s+1}}-
 \f{ax+by}{3^s})\o}2^{q_1}\Cs{-2\f{ax+by}{3^s}+(\f{Py}{3^{s+1}}-
 \f{ax+by}{3^s})\o}y^2
 \\&\q\times\Cs{-6b+(P-3b)\o^2}{a}
 \\&=\Cs{\o}{ap}^2\Cs{\o}2^{q_1}\Cs{-6(ax+by)+(Py-3(ax+by)\o}y^2
 \Cs{P-3b-6b\o}{a}
 \\&=\Cs{\o}{ap}^2\cdot \o^{q_1}\Cs{-3ax(2+\o)}y^2\Cs{P-3b-6b\o}{a}
 \\&=\Cs{\o}{ap}^2\cdot \o^{q_1}\Cs{\o(1+2\o)}{y^2}\Cs{P-3b(1+2\o)}a
 \\&=\Cs{\o}{ap}^2\cdot \o^{q_1}\cdot \o^{\f{1-y^2}3}\Cs{P-3b(1+2\o)}a
 =\Cs{P-3b(1+2\o)}a.\end{align*}
\par Summarizing the above proves the theorem.
\vskip0.2cm
\par{\bf Remark 2.1} By Case 1 in the proof of Theorem 2.1, if $P\e \pm 3\mod 9$ and $y\e \pm 9\mod {27}$, or if $P\e \pm 9\mod {27}$ and $y\e \pm 3\mod 9$, we have
$\cs{Py+3(ax+by)(1+2\o)}{ax^2+2bxy+cy^2}\not=\cs{P-3b(1+2\o)}a.$

\section*{3. The cubic character on $H(4D)$}

\par{\bf Theorem 3.1} {\sl Let $P,D\in\Bbb Z$ with $3\nmid P$ or $27\mid P$. Suppose $Q=P^2+27D=2^{q_1}q_0^3$ with $q_0\in\Bbb Z$ and $q_1\in\{0,1,2\}$. Assume that $DQ\not=0$. For $[a,2b,c]\in H(4D)$ with $(a,3DQ)=1$ define
$\chi([a,2b,c])=\cs{P-3b(1+2\o)}a$. Then $\chi$ is
 a character on $H(4D)$ and so
 $$G(P,D)=\Big\{[a,2b,c]\in H(4D)\Bigm| (a,3DQ)=1, \Cs{P-3b(1+2\o)}a=1\Big\}$$
is a subgroup in $H(4D)$, where $a>0$ when $D<0$.
Moreover, if the congruence $Dx^3-Qx-2Q\e 0\mod {p_0}$ is insolvable for some odd prime $p_0$ with $p_0\nmid 3DPQ$, then $G(P,D)$ is a subgroup of index $3$ in $H(4D)$.}
\v2
\par{\it Proof.} We first illustrate that $\chi$ is well defined.
For any given positive integer $M$ and a class $K\in H(4D)$, from [2,16] we may assume that $K=[a,2b,c]$ with $(a,M)=1$.
Suppose that $(a,2b,c)$ and $(a',2b',c')$ are equivalent with  $b^2-ac={b'}^2-a'c'=D$. Then
there are integers
$\al,\be,\ga,\de$ such that $\al\de-\be\ga=1$ and
$$\align a(\al X&+\be Y)^2+2b(\al X+\be Y)(\ga X+\de Y)+c(\ga X+\de Y)^2
\\&=a'X^2+2b'XY+c'Y^2.\endalign$$
Hence
\begin{align*} &a'=a\al^2+2b\al\ga+c\ga^2,\q
b'=a\al\be+b(\al\de+\be\ga)+c\ga\de,
\\&c'=a\be^2+2b\be\de+c\de^2.\end{align*}
 Set $a^*=a/(a,\ga),\ c^*=(a,\ga)c,\ x=\a$
and $y=\ga/(a,\ga)$. Then
$$a^*x^2+2bxy+c^*y^2=\f{a\al^2+2b\al\ga+c\ga^2}{(a,\ga)}=\f{a'}{(a,\ga)}$$
and
$$\align b'\ga&=a\al\be\ga+b(\al\ga\de+\be\ga^2)+c\ga^2\de
\\&\e a\al\be\ga+b(\al\ga\de+\be\ga^2)-\de(a\al^2+2b\al\ga)
\\&=(\be\ga-\al\de)(a\al+b\ga)=-a\a-b\ga\mod{|a'|}.\endalign$$
Hence $$\f{b'\ga}{(a,\ga)}\e
-\f{a\a}{(a,\ga)}-b\f{\ga}{(a,\ga)}
=-a^*x-by\mod{\f{|a'|}{(a,\ga)}}.$$
Since $(aa',3DQ)=1$ we have $(a^*(a^*x^2+2bxy+c^*y^2),3DQ)
=(\f a{(a,\ga)}\cdot\f{a'}{(a,\ga)},3DQ)=1$. Also,
$(\f a{(a,\ga)},\f{\ga}{(a,\ga)})=1$ and
$$\Big(\f {a'}{(a,\ga)},\f{\ga}{(a,\ga)}\Big)
=\Big(\f a{(a,\ga)}\a^2+2b\a\f{\ga}{(a,\ga)}+c\f{{\ga}^2}{(a,\ga)},\f{\ga}{(a,\ga)}\Big)
=\Big(\f a{(a,\ga)}\a^2,\f{\ga}{(a,\ga)}\Big)=1.$$
Thus, $(a^*(a^*x^2+2bxy+c^*y^2),y)=1$. Now, from the above and Theorem 2.1 we deduce that
\begin{align*}\Cs{P-3b'(1+2\o)}{a'/(a,\ga)}&
=\Cs{P\f{\ga}{(a,\ga)}-3b'\f {\ga}{(a,\ga)}(1+2\o)}{a'/(a,\ga)}
\\&=\Cs{Py+3(a^*x+by)(1+2\o)}{a^*x^2+2bxy+c^*y^2}
=\Cs{P-3b(1+2\o)}{a^*}.
\end{align*}
On the other hand, since $b'\e b\a\de=b(1+\be\ga)\e b\mod {(a,\ga)}$ we have
$\cs{P-3b'(1+2\o)}{(a,\ga)}=\cs{P-3b(1+2\o)}{(a,\ga)}$ and so
\begin{align*}\Cs{P-3b'(1+2\o)}{a'}
&=\Cs{P-3b'(1+2\o)}{a'/(a,\ga)}
\Cs{P-3b'(1+2\o)}{(a,\ga)}
\\&=\Cs{P-3b(1+2\o)}{a/(a,\ga)}\Cs{P-3b(1+2\o)}{(a,\ga)}
\\&=\Cs{P-3b(1+2\o)}{a}.\end{align*}
Hence $\chi$ is well defined.
\par For $[A_1,2B_1,C_1],[A_2,2B_2,C_2]\in H(4D)$ with $(A_1A_2, 3DQ)=1$ let
$$t=\t{\rm gcd}(A_1,A_2,B_1+B_2)=A_1u+A_2v+(B_1+B_2)w\ (u,v,w\in\Bbb Z),$$ and let
$$A_3=\f{A_1A_2}{t^2}, \q B_3=B_2+\f{A_2}t((B_1-B_2)v-C_2w)\qtq{and}C_3=\f{B_3^2-D}{A_3}.$$
From [1,\;p.246], $[A_1,2B_1,C_1][A_2,2B_2,C_2]=[A_3,2B_3,C_3]$. Clearly
 $B_3\e B_2\mod {\f{|A_2|}t}$. Since $C_2=\f{B_2^2-D}{A_2}=\f{B_2^2-B_1^2+A_1C_1}{A_2}$, we also have
 \begin{align*}B_3&=B_2+\f{A_2}t\Big((B_1-B_2)v-\f{B_2^2-B_1^2+A_1C_1}{A_2}w\Big)
 \\&=B_2+(B_1-B_2)\Big(\f{A_2}tv+\f{B_1+B_2}tw\Big)-\f{A_1}tC_1w
 \\&=B_2+(B_1-B_2)\Big(1-\f{A_1}tu\Big)-\f{A_1}tC_1w
 \\&\e B_2+(B_1-B_2)=B_1\mod {\f{|A_1|}t}.
 \end{align*}
It is clear that for $i=1,2$,
\begin{align*}N(P-3B_i(1+2\o))&=(P-3B_i)^2+(P-3B_i)\cdot 6B_i+(-6B_i)^2
=P^2+27B_i^2
\\&=P^2+27(D+A_iC_i)=Q+27A_iC_i.\end{align*}
Recall that $(A_i,Q)=1$. We then have $(A_i,N(P-3B_i(1+2\o)))=1$
and so $\cs{P-3B_i(1+2\o)}{A_i}\not=0$ for $i=1,2$.
In view of $B_1+B_2\e 0\mod t$,
$$\Cs{P-3B_1(1+2\o)}t\Cs{P-3B_2(1+2\o)}t
=\Cs{P^2+3B_1\cdot 3B_2\cdot(-3)}t=1.$$
Now, from the above we derive that
\begin{align*}\Cs{P-3B_3(1+2\o)}{A_3}
&=\Cs{P-3B_3(1+2\o)}{A_1/t}\Cs{P-3B_3(1+2\o)}{A_2/t}
\\&=\Cs{P-3B_1(1+2\o)}{A_1/t}\Cs{P-3B_2(1+2\o)}{A_2/t}
\\&=\Cs{P-3B_1(1+2\o)}{A_1}\Cs{P-3B_2(1+2\o)}{A_2}
\\&\q\times\Cs{P-3B_1(1+2\o)}t^{-1}\Cs{P-3B_2(1+2\o)}t^{-1}
\\&=\Cs{P-3B_1(1+2\o)}{A_1}\Cs{P-3B_2(1+2\o)}{A_2}.
\end{align*}
Also, $\chi([1,0,-D])=\cs{P}1=1$. Hence $\chi$ is a character on $H(4D)$ and the kernel $G(P,D)$ is a subgroup of $H(4D)$.
\par Now assume that the congruence $Dx^3-Qx-2Q\e 0\mod {p_0}$ is insolvable for some odd prime $p_0$ with $p_0\nmid 3DPQ$. Since the discriminant of $x^3-\f QDx-\f{2Q}D$ is $\f{4P^2Q^2}{D^3}$, it follows from (1.2) that $\sls D{p_0}=1$. Suppose $b^2\e D\mod {p_0}$ $(b\in\Bbb Z)$. Then $x^3-9(\f{P^2}{9b^2}+3)x-18(\f{P^2}{9b^2}+3)\e 0\mod {p_0}$ is insolvable. Applying [12, Lemma 4.1] gives $\cs{P-3b(1+2\o)}{p_0}\not=1$. Set $c=(b^2-D)/p_0$.
Then $[p_0,2b,c]\in H(4D)$ and $(p_0,3DQ)=1$. Since $\chi([p_0,2b,c])=\o$ or $\o^2$, from the above we see that
$$\{\chi([p_0,2b,c]),\ \chi([p_0,2b,c]^2),\chi([p_0,2b,c]^3)\}=\{\o,\o^2,1\}.$$
Therefore $\chi$ is a surjective homomorphism from $H(4D)$ to $\{1,\o,\o^2\}$ and so the kernel
$G(P,D)$ is a subgroup of index $3$ in $H(4D)$. The proof is now complete.

\v2
\par{\bf Corollary 3.1} {\sl Let $a_1,a_2,a_3\in\Bbb Z$ and $P_0,Q_0,D_0,P_1,Q_1,D_1$ be as in (1.1) and (1.3).  For $[a,2b,c]\in H(4D_1)$ with $(a,3D_1Q_1)=1$ define
$\chi([a,2b,c])=\cs{P_1-3b(1+2\o)}a$. If there is a prime $p_0$ such that $p_0\nmid 3D_1Q_1$ and the congruence $D_1x^3-Q_1x-2Q_1\e 0\mod {p_0}$ is insolvable, then $\chi$ is a surjective homomorphism from $H(4D_1)$ to $\{1,\o,\o^2\}$ and so the kernel
$$G(P_1,D_1)
=\Big\{[a,2b,c]\in H(4D_1)\Bigm| (a,3D_1Q_1)=1, \Cs{P_1-3b(1+2\o)}a=1\Big\}$$
is a subgroup of index $3$ in $H(4D_1)$, where $a>0$ when $D_1<0$.}
\v2
\par{\it Proof.} If $3\mid a_1$, then $P_0=-2a_1^3+9a_1a_2-27a_3\e 0\mod {27}$ and so $27\mid P_1$. If $3\nmid a_1$, then $3\nmid P_0$ and so $3\nmid P_1$. Now the result follows from Theorem 3.1 immediately.

\section*{4. Criteria for $N_p(x^3+a_1x^2+a_2x+a_3)=3$}
\par{\bf Theorem 4.1} {\sl Let $a_1,a_2,a_3\in\Bbb Z$ and $P_0,Q_0,D_0,P_1,Q_1,D_1$ be given in (1.1) and (1.3). Let
$$G(a_1,a_2,a_3)
=\Big\{[a,2b,c]\in H(4D_1)\Bigm| (a,3D_1Q_1)=1, \Cs{P_1-3b(1+2\o)}a=1\Big\},$$
where $a>0$ if $D_1<0$.
Suppose that $p>3$ is a prime such that $p\nmid D_0Q_0$. Then the following statements are equivalent}:
\par (i) $x^3+a_1x^2+a_2x+a_3\e 0\mod p$ {\it has three solutions};
\par (ii) {\it $p$ is represented by some class in $G(a_1,a_2,a_3)$};
\par (iii) $s_{p+1}(a_1,a_2,a_3)\e a_1^2-2a_2\mod p$;
\par (iv) $V_{\f{p-\sls p3}3}(P_0,Q_0)\e 2(a_1^2-3a_2)^{\f{1-\sls p3}2}
\mod p$.
\newline {\t Moreover, if $p\nmid P_0$, then the above statements are equivalent to}
\par (v) $u_{p-2}(a_1,a_2,a_3)\e 0\mod p$;
\par (vi) $U_{\f{p-\sls p3}3}(P_0,Q_0)\e 0\mod p$;
 \newline \t{if $p\nmid P_0(P_0^2-3Q_0)$, the above statements are equivalent to}
 \par (vii) $U_{2[\f p3]+1}(P_0,Q_0)\e (-Q_0)^{[\f p3]}\mod p$;
\par (viii) $\sum_{k=1}^{[p/3]}\b{3k}k\sls{P_0^2}{27Q_0}^k\e 0\mod p$.
\v2
\par{\it Proof.} If $p=ax^2+2bxy+cy^2$ with $a,b,c,x,y\in\Bbb Z$
  and $b^2-ac=D_1$, then $ap=(ax+by)^2-D_1y^2$. It is clear that $(x,y)=1$ since $(x,y)\mid p$ and $p^2\nmid p$.
  For any given positive integer $M$ and a class $K\in H(4D_1)$, from [2,16] we may assume that $K=[a,2b,c]$ with $(a,M)=1$. Hence we may assume $K=[a,2b,c]$ with $(a,3pD_1Q_1)=1$. For $p\nmid a$ we have
  $p\nmid y$ since $p\mid y$ implies that $p\mid ax$ and so $p\mid a$. Hence $(ap,y)=1$ and $(\f{ax+by}y)^2\e D_1\mod p$. Since $3\nmid P_1$ or $27\mid P_1$ according as $3\nmid a_1$ or $3\mid a_1$, from [13, Lemma 4.2] and Theorem 2.1 we deduce that
\begin{align*} &\t{$x^3+a_1x^2+a_2x+a_3\e 0\mod p$ has three solutions}
\\&\Leftrightarrow \t{there is an integer $d_0$ such that $d_0^2\e D_0\mod p$ and $\cs{P_0+3d_0(1+2\o)}p=1$}
\\&\Leftrightarrow \t{there is an integer $d_1$ such that $d_1^2\e D_1\mod p$ and $\cs{P_1+3d_1(1+2\o)}p=1$}
\\&\Leftrightarrow \t{$p=ax^2+2bxy+cy^2$ with $a,b,c,x,y\in\Bbb Z$,
  $(a,3pD_1Q_1)=1$,\ \t{\rm gcd}$(a,2b,c)=1$, }
 \\&\qq b^2-ac=D_1\q \t{and}\q
  \Cs{P_1y+3(ax+by)(1+2\o)}p=1
 \\&\Leftrightarrow \t{$p=ax^2+2bxy+cy^2$ with $a,b,c,x,y\in\Bbb Z$,   $(a,3pD_1Q_1)=1$,}
 \\&\qq \t{$\t{\rm gcd}(a,2b,c)=1$, $b^2-ac=D_1$ and $\Cs{P_1-3b(1+2\o)}a=1$.}
 \end{align*}
 If $(a,2b,c)$ and $(a',2b',c')$ are equivalent forms with discriminant $4D_1$ and $(a,3pD_1Q_1)=(a',3D_1Q_1)=1$, from the proof of Theorem 3.1 we know that $\cs{P_1-3b(1+2\o)}a=\cs{P_1-3b'(1+2\o)}{a'}$. Thus, from the above we deduce that
  \begin{align*}
 &\t{$x^3+a_1x^2+a_2x+a_3\e 0\mod p$ has three solutions}
 \\&\Leftrightarrow \t{$p=ax^2+2bxy+cy^2$ with $a,b,c,x,y\in\Bbb Z$,   $(a,3D_1Q_1)=1$,}
 \\&\qq \t{$\t{\rm gcd}(a,2b,c)=1$, $b^2-ac=D_1$ and $\Cs{P_1-3b(1+2\o)}a=1$.}
 \\&\Leftrightarrow \t{$p$ is represented by a class in $G(a_1,a_2,a_3)$.}
\end{align*}
This shows that (i) is equivalent to (ii). By [13, Theorem 4.1], (i) is equivalent to (iii).  From [13, Lemma 3.1], (iii) is equivalent to (iv).
\par Now assume that $p\nmid P_0$. From [13, Theorem 4.3],
(i) is equivalent to (v). By [13, Theorem 3.2(i)], (v) is equivalent to (vi). From [15, Lemmas 2.2 and 2.4], (vi),(vii) and (viii) are equivalent when $p\nmid P_0^2-3Q_0$. Thus the theorem is proved.
\vskip0.2cm
\par{\bf Theorem 4.2} {\sl Let $a_2,a_3\in\Bbb Z$ and let $p>3$ be a prime such that $p\nmid a_2a_3(4a_2^3+27a_3^2)$. Then
$p\mid u_p(0,a_2,a_3)$ if and only if $p$ is represented by a class in
\begin{align*}G(a_2,a_3)&=\Big\{[a,2b,c]\in H(-2^{1-(-1)^{a_3}}(4a_2^3+27a_3^2))\Bigm| \Big(a,3a_2\f{4a_2^3+27a_3^2}{2^{1+(-1)^{a_3}}}\Big)=1, \\&\qq \Cs{b-\f {6a_3}{3+(-1)^{a_3}}(1+2\o)}a=1\Big\},\end{align*}
where $a>0$ for $4a_2^3+27a_3^2>0$.}
\par{\it Proof.}
By [13, Theorem 3.3], if $\ls{-4a_2^3-27a_3^2}p=-1$, then
$$-(4a_2^3+27a_3^2)u_p(0,a_2,a_3)\e -9a_2a_3x^2-2a_2^3x-6a_2^2a_3\mod p,$$
where $x$ is an integer satisfying $x^3+a_2x+a_3\e 0\mod p$.
Thus,
\begin{align*}
p\mid u_p(0,a_2,a_3)&\Leftrightarrow 9a_3x^2+2a_2^2x+6a_2a_3\e 0\mod p
\\&\Leftrightarrow 9a_3x^3+2a_2^2x^2+6a_2a_3x\e 0\mod p
\\&\Leftrightarrow -9a_3(a_2x+a_3)+2a_2^2x^2+6a_2a_3x\e 0\mod p
\\&\Leftrightarrow (2a_2x+3a_3)(a_2x-3a_3)\e 0\mod p
\\&\Leftrightarrow x\e -\f{3a_3}{2a_2}\ \t{or}\ \f{3a_3}{a_2}\mod p.
\end{align*}
Note that
\begin{align*} &\Big(-\f{3a_3}{2a_2}\Big)^3+a_2\Big(-\f{3a_3}{2a_2}\Big)+a_3
=-\f{a_3(4a_2^3+27a_3^2)}{8a_2^3}\not\e 0\mod p,
\\&\Big(\f{3a_3}{a_2}\Big)^3+a_2\cdot \f{3a_3}{a_2}+a_3
=\f{a_3(4a_2^3+27a_3^2)}{a_2^3}\not\e 0\mod p.
\end{align*}
We then have $p\nmid u_p(0,a_2,a_3)$ for $\ls{-4a_2^3-27a_3^2}p=-1$.
By (1.2), $N_p(x^3+a_2x+a_3)=1\Leftrightarrow \ls{-4a_2^3-27a_3^2}p=-1$.
Hence, $N_p(x^3+a_2x+a_3)=1$ implies that $p\nmid u_p(0,a_2,a_3)$.
\par If $\ls{-4a_2^3-27a_3^2}p=1$ and $d^2\e -4a_2^3-27a_3^2\mod p$,
it follows from [13, Theorem 3.3 and Lemma 4.2] that
$$p\mid u_p(0,a_2,a_3)\Leftrightarrow \Cs{-27a_3+3d(1+2\o)}p=1
\Leftrightarrow N_p(x^3+a_2x+a_3)=3.$$
\par Now, from the above and Theorem 4.1 we deduce that
\begin{align*}p\mid u_p(0,a_2,a_3)
&\Leftrightarrow N_p(x^3+a_2x+a_3)=3
\\&\Leftrightarrow \t{$p$ is represented by a class in $G(0,a_2,a_3)$}.
\end{align*}
Set $a_3'=\f{2a_3}{3+(-1)^{a_3}}$. Then
 $$\Cs{-27a_3'-3b(1+2\o)}a=\Cs{-3(1+2\o)}a\cs{b-3a_3'(1+2\o)}a
 =\Cs{b-3a_3'(1+2\o)}a.$$
 Thus
 $G(0,a_2,a_3)=G(a_2,a_3)$ and the theorem is proved.
\v2
\par{\bf Corollary 4.1} {\it Let $a_1,a_2,a_3\in\Bbb Z$ and let $p>3$ be a prime such that $p\nmid D_0P_0Q_0$, where $P_0,Q_0$ and $D_0$ are given in (1.1). Then}
$$N_p(x^3+a_1x^2+a_2x+a_3)=3\Leftrightarrow p\mid u_p(0,-3(a_1^2-3a_2), 2a_1^3-9a_1a_2+27a_3).$$
\par{\it Proof.} Set $P_0=-2a_1^3+9a_1a_2-27a_3$. Observe that
$$27(x^3+a_1x^2+a_2x+a_3)=(3x+a_1)^3-3(a_1^2-3a_2)(3x+a_1)-P_0.$$
We have
$$N_p(x^3+a_1x^2+a_2x+a_3)=N_p(x^3-3(a_1^2-3a_2)x-P_0).$$
By the proof of Theorem 4.2, $N_p(x^3-3(a_1^2-3a_2)x-P_0)=3$ if and only if $p\mid u_p(0,-3(a_1^2-3a_2), -P_0)$. Thus the result follows.
\v2
\par{\bf Example 4.1}  Let $p$ be a prime with $p\not=2,3,5,7,13$. Then
\begin{align*}&x^3-x+2\e 0\mod p\q\t{has three solutions}
\\&\Leftrightarrow p=x^2+26y^2\ \t{or}\ 2x^2+13y^2\ (x,y\in\Bbb Z)
\\&\Leftrightarrow s_{p+1}(0,-1,2)\e 2\mod p
\Lra p\mid u_{p-2}(0,-1,2)
\Lra p\mid u_p(0,-1,2)
\\&\Lra p\mid U_{\f{p-\sls p3}3}(18,3)
\Lra \sum_{k=1}^{[p/3]}\b{3k}k4^k\e 0\mod p.
\end{align*}
\par{\it Proof.} Set $a_1=0,\ a_2=-1$ and $a_3=2$. Then $P_0=-54$, $Q_0=27$, $D_0=-104$, $P_1=-27,\ Q_1=27$ and $D_1=-26$ by (1.1) and (1.3).  Calculating reduced forms with discriminant $-104$, we see that
$$H(-104)=\{[1,0,26],\ [2,0,13], \ [3,2,9],\ [3,-2,9],\ [5,4,6],\ [5,-4,6]\}.$$
It is well known that (see [2]) $(a,b,c)\sim (c,-b,a)$ and
$(a,b,c)\sim (a,2ak+b,ak^2+bk+c)$ for $k\in\Bbb Z$. Thus,
$[2,0,13]=[2,12,31]=[31,-12,2]$ and $[3,2,9]=[3,-10,17]=[17,10,3]$.
Using (2.1)-(2.6) we find that
\begin{align*}&\Cs{-27-3\cdot 0(1+2\o)}{1}=1,\ \Cs{-27-3(-6)(1+2\o)}{31}=1,\\& \Cs{-27-3\cdot 5(1+2\o)}{17}=\o,
\ \Cs{-27-3\cdot 2(1+2\o)}{5}=\o.\end{align*}
Thus,
$G(0,-1,2)=\{[1,0,26],[31,-12,2]\}=\{[1,0,26],[2,0,13]\}.$
Now applying Theorems 4.1-4.2 and the fact that $U_n(-54,27)=(-3)^{n-1}U_n(18,3)$ yields the result.
\v2
\par{\bf Example 4.2}  Let $p$ be a prime with $p\not=2,3,7,11,79,3119$. Then
\begin{align*}&x^3+2x^2+x+3\e 0\mod p\q\t{has three solutions}
\\&\Leftrightarrow p=x^2+231y^2,\ 3x^2+77y^2,\ 7x^2+33y^2\
\t{or}\ 11x^2+21y^2\ (x,y\in\Bbb Z)
\\&\Leftrightarrow s_{p+1}(2,1,3)\e 2\mod p
\Lra p\mid u_{p-2}(2,1,3)
\Lra p\mid u_p(0,-3,79)
\\&\Lra p\mid U_{\f{p-\sls p3}3}(-79,1)
\Lra \sum_{k=1}^{[p/3]}\b{3k}k\Ls{6241}{27}^k\e 0\mod p.
\end{align*}
\par{\it Proof.} Set $a_1=2,\ a_2=1$ and $a_3=3$. Then $P_0=-79$, $Q_0=1$,
   $D_0=-231$, $P_1=-79,\ Q_1=4$ and $D_1=-231$ by (1.1) and (1.3).
   Computing reduced forms with discriminant $-924$, we get
\begin{align*}H(-924)=\{&[1,0,231],\ [3,0,77], \ [7,0,33],\ [11,0,21],\ [8,2,29],\ [8,-2,29],
\\&[5,4,47],\ [5,-4,47],\ [15,6,16],\ [15,-6,16],\ [13,8,19],\ [13,-8,19]\}.
\end{align*}
Since $[3,0,77]=[3,-12,89]=[89,12,3]$, $[7,0,33]=[7,-28,61]=[61,28,7]$,
$[11,0,21]=[11,-88,197]=[197,88,11]$ and
\begin{align*}&\Cs{-79-3\cdot 0(1+2\o)}{1}=1,\ \Cs{-79-3\cdot 6(1+2\o)}{89}=1,\\& \Cs{-79-3\cdot 14(1+2\o)}{61}=1,
\ \Cs{-79-3\cdot 44(1+2\o)}{197}=1,\end{align*}
we see that
\begin{align*}G(2,1,3)&=\{[1,0,231], \ [89,12,3],\ [61,28,7],\ [197,88,11]\}
\\&=\{[1,0,231],\ [3,0,77], \ [7,0,33],\ [11,0,21]\}.
\end{align*}
Now applying Theorem 4.1 and Corollary 4.1 gives the result.

\section*{5. The criterion for $Q(P+\sqrt{P^2-Q})$ to be a cubic residue $\mod p$}
\par\q\ {\bf Lemma 5.1} {\sl Suppose that $p$ is a prime of the form $3k+1$, $P,Q,D\in\Bbb Z$, $p\nmid DQ$, $P^2+27D=Q$, $\ls DP=1$ and $d^2\e D\mod p$. Then $Q(P+3\sqrt{-3D})$ is a cubic residue $\mod p$ if and only if
$\cs{P+3d(1+2\o)}p=1$.}
\v2
\par{\it Proof.} Let $t\in\Bbb Z$ be such that $t^2\e -3\mod p$. From [12, Theorem 2.2],
\begin{align*}\Cs{P+3d(1+2\o)}p=1&
\Lra \Big(\f{\f P{3d}-t}{\f P{3d}+t}\Big)^{\f{p-1}3}\e 1\mod p
\\&\Lra\f{P-3dt}{P+3dt}\q\t{is a cubic residue $\mod p$}
\\&\Lra\f{P-3\sqrt{-3D}}{P+3\sqrt{-3D}}\q\t{is a cubic residue $\mod p$}
\\&\Lra \f Q{(P+3\sqrt{-3D})^2}\q\t{is a cubic residue $\mod p$}
\\&\Lra Q(P+3\sqrt{-3D})\q\t{is a cubic residue $\mod p$.}
\end{align*}
\par{\bf Theorem 5.1} {\sl Let $P\in\Bbb Z$ with $3\nmid P$ or $27\mid P$, and let $D,Q\in\Bbb Z$ be such that $Q=P^2+27D=2^{q_1}q_0^3$ with
 $q_1\in\{0,1,2\}$ and $q_2\in\Bbb Z$.
Suppose that $p$ is a prime of the form $3k+1$, $p\nmid Q$ and $\ls {P^2-Q}p=1$. Then $Q(P+\sqrt{P^2-Q})$ is a cubic residue $\mod p$ if and only if $p$ is represented by some class in
$G(P,D)$, where $G(P,D)$ is given in Theorem 3.1.}
\v2
\par{\it Proof.} If $p=ax^2+2bxy+cy^2$ with $a,b,c,x,y\in\Bbb Z$, $(a,3pDQ)=1$ and $b^2-ac=D$, then clearly $p\nmid y$ and $\ls{ax+by}y^2\e D\mod p$. Appealing to Lemma 5.1, Theorem 2.1 and the argument in the proof of Theorem 4.1 we deduce that
\begin{align*} &Q(P+\sqrt{P^2-Q})\q\t{is a cubic residue $\mod p$}
\\&\Lra Q(P+3\sqrt{-3D})\q\t{is a cubic residue $\mod p$}
\\&\Leftrightarrow \t{for $d\in\Bbb Z$ with $d^2\e D\mod p$ we have $\cs{P+3d(1+2\o)}p=1$}
\\&\Leftrightarrow \t{$p=ax^2+2bxy+cy^2$ with $a,b,c,x,y\in\Bbb Z$,
  $(a,3pDQ)=1$,\ \t{\rm gcd}$(a,2b,c)=1$, }
 \\&\qq b^2-ac=D\q \t{and}\q\Cs{Py+3(ax+by)(1+2\o)}p=1
 \\&\Leftrightarrow \t{$p=ax^2+2bxy+cy^2$ with $a,b,c,x,y\in\Bbb Z$,   $(a,3pDQ)=1$,}
 \\&\qq \t{$\t{\rm gcd}(a,2b,c)=1$, $b^2-ac=D$ and $\Cs{P-3b(1+2\o)}a=1$}
 \\&\Leftrightarrow \t{$p=ax^2+2bxy+cy^2$ with $a,b,c,x,y\in\Bbb Z$,   $(a,3DQ)=1$,}
 \\&\qq \t{$\t{\rm gcd}(a,2b,c)=1$, $b^2-ac=D$ and $\Cs{P-3b(1+2\o)}a=1$}
 \\&\Leftrightarrow \t{$p$ is represented by a class in $G(P,D)$.}
\end{align*}

\par{\bf Example 5.1} {\sl Let $p$ be a prime such that $p\e 1\mod 3$ and $\ls{87}p=1$. Then
$28+3\sqrt{87}$ is a cubic residue $\mod p$ if and only if $p=x^2+29y^2$ or $2x^2+2xy+15y^2$ for some $x,y\in\Bbb Z$.}
\par{\it Proof.} Taking $P=28,\ Q=1,\ D=-29$ in Theorem 5.1 and noting that $H(-116)=\{[1,0,29],[2,2,15],[3,2,10],[3,-2,10],[5,2,6],[5,-2,6]\}$ yields the result.
\v2
\par{\bf Remark 5.1} Putting $d=87$ in [4, Theorem 6.1] yields that $28+3\sqrt{87}$ is a cubic residue $\mod p$ for primes $p=x^2+29y^2\e 1\mod 3$.
\v2
\par{\bf Example 5.2} {\sl Let $p$ be a prime such that $p\e 1\mod 3$ and $\ls{114}p=1$. Then
$64+6\sqrt{114}$ is a cubic residue $\mod p$ if and only if $p=x^2+38y^2$ or $2x^2+19y^2$ for some $x,y\in\Bbb Z$.}
\par{\it Proof.} Taking $P=32,\ Q=-2,\ D=-38$ in Theorem 5.1 and noting that $H(-152)=\{[1,0,38],[2,0,19],[3,2,13],[3,-2,13],[6,4,7],[6,-4,7]\}$ yields the result.
\v2
\par{\bf Example 5.3} {\sl Let $p$ be a prime such that $p\e 1\mod 3$ and $\ls{78}p=1$. Then
$27+3\sqrt{78}$ is a cubic residue $\mod p$ if and only if $p=x^2+26y^2$ or $2x^2+13y^2$ for some $x,y\in\Bbb Z$.}
\par{\it Proof.} Taking $P=Q=27,\ D=-26$ in Theorem 5.1 and noting that $H(-104)=\{[1,0,26],[2,0,13],[3,2,9],[3,-2,9],[5,4,6],[5,-4,6]\}$ yields the result.
\v2
\par{\bf Acknowledgements.} The author was supported by the
National Natural Science Foundation of China (Grant No. 12271200).


\begin{thebibliography} {99}

\bibitem [1] {}  H.
Cohen, {\it A Course in Computational Algebraic Number Theory},
      Graduate Texts in Mathematics 138, Springer, Berlin, New York, 1993.

\bibitem [2] {} D.A. Cox, {\it Primes of the Form $x^2+ny^2$:
Fermat, Class Field Theory, and Complex Multiplication},
Wiley, Inc., New York, Chichester, 1989.

\bibitem [3] {} L.E. Dickson, {\it Criteria for the irreducibility
of functions in a finite field}, Bull. Amer. Math. Soc.
13(1906), 1-8.

\bibitem [4] {}
R. Evans, F. Lemmermeyer, Z.H. Sun and M.V. Veen, {\it Ring class fields and a result of Hasse}, J. Number Theory 266(2025), 33-61.

\bibitem [5] {} T. Evink and P.A. Helminck, {\it Tribonacci numbers and primes of the form $p=x^2+11y^2$}, Math.
Slovaca {\bf 69}(2019), 521-532.

\bibitem [6] {} A. Faisant, {\it On the Padovan sequence}, preprint (2019), arXiv:1905.07702.

\bibitem [7] {} K. Ireland and
M. Rosen, {\it A Classical Introduction to Modern Number Theory}, second edition, Springer, New York, 1990.

\bibitem [8] {} P. Moree and A. Noubissie, {\it Higher reciprocity laws and ternary linear recurrence sequences}, preprint (2022), arXiv:2205.06685.

\bibitem [9] {} T. Skolem, {\it On a certain connection between
the discriminant of a polynomial and the number of its irreducible
factors mod p}, Norsk Mat. Tidsskr. 34(1952), 81-85.

\bibitem [10] {} B.K. Spearman and K.S. Williams,
{\it The cubic congruence $x^3+Ax^2+Bx+C\e 0\mod p$ and binary
quadratic forms}, Proc. London Math. Soc. 46(1992), 397-410.

\bibitem [11] {} B.K. Spearman and K.S. Williams, {\it The cubic congruence
$x^3+Ax^2+Bx+C\e 0\mod p$ and binary quadratic forms II}, J.
London Math. Soc. 64(2001), 273-274.

\bibitem [12] {} Z.H. Sun,
{\it On the theory of cubic residues and nonresidues}, Acta
Arith. 84(1998), 291-335.

 \bibitem [13] {} Z.H. Sun,
{\it Cubic and quartic congruences modulo a prime}, J.
Number Theory 102(2003), 41-89.


     \bibitem [14] {} Z.H. Sun,  {\it Cubic residues and binary quadratic forms},
 J. Number Theory 124(2007), 62-104.

 \bibitem [15] {} Z.H. Sun, {\it Cubic congruences and sums involving $\b{3k}k$}, Int. J. Number Theory 12(2016), 143-164.

 \bibitem [16] {} Z.H. Sun and K. S. Williams,
{\it On the number of representations of $n$ by $ax^2+
bxy+cy^2$}, Acta Arith.  122(2006), 101-171.

\end{thebibliography}
\end{document}